\documentclass[a4paper, 11pt]{amsart}
\usepackage[margin=2.3cm]{geometry}
\usepackage{tikz}
\usepackage{wrapfig}
\usetikzlibrary{arrows} 
\usepackage[margin=0pt]{subfig}
\theoremstyle{plain}
\newtheorem{theorem}{Theorem}[section]

\theoremstyle{definition}

\newcommand{\R}{{\mathbb{R}}}

\newcommand{\N}{{\mathbb{N}}}

\numberwithin{equation}{section}

\begin{document}

%\volumetitle{ICM 2022} % Don't alter this line.

%------
% Insert the title of your paper and (if necessary)
% a short title for the running head.
%------

\title[Dynamics and 'arithmetics' of surface flows]{Dynamics and `arithmetics' of higher genus surface flows}

\author[C.\ Ulcigrai]{Corinna Ulcigrai}
\address{Institut f\"ur Mathematik, Universit\"at Z\"urich, Winterthurerstrasse 190,
CH-8057 Z\"urich, Switzerland}
\email{corinna.ulcigrai@math.uzh.ch}\date{}

\begin{abstract}
We survey some recent advances in the study of (area-preserving) flows on surfaces, in particular on the typical dynamical, ergodic and spectral properties of smooth area-preserving (or \emph{locally Hamiltonian}) flows, as well as recent breakthroughs on  \emph{linearization} and \emph{rigidity} questions in higher genus. We focus in particular on the \emph{Diophantine-like conditions} which are required to prove such results, which can be thought of as a generalization of   \emph{arithmetic conditions} for flows on tori and circle diffeomorphisms. We will explain how 
these conditions on higher genus flows and their Poincar{\'e} sections (namely generalized interval exchange maps) can be imposed by controlling a  renormalization dynamics, but are of more subtle nature than in genus one since they often exploit features which originate from  the non-uniform hyperbolicity of the renormalization. % nature.
%  a role analogous to that of \emph{arithmetic conditions} for flows on tori and circle diffeomorphisms. As classical arithmetic conditions, Diophantine-like conditions in higher genus  but their nature is more suble. 
%Diophantine and more in general conditions in genus one can be expressed in terms of continued fraction entries, \emph{Diophantine-like conditions} in higher genus are expressed using renormalization tools and are of more subtle nature, as we will explain in this survey. 
\end{abstract}

\maketitle

%------
% INSERT THE BODY OF THE PAPER HERE (except
% acknowledgments, funding info and bibliography)
%------

\section{Introduction.}\label{sec:intro}
Flows on surfaces are one of the most basic and most fundamental examples of dynamical systems. First of all, 
%both for their \emph{low dimension} and for their appearence in problems of physical origin. 
%Smooth surface flows are fundamental in dynamics is because 
they are  among the lowest possible dimensional smooth systems; % (on compact manifolds of lower dimension, the  other fundamental class of smooth dynamical systems are circle diffeomorphisms, whose rich theory is a cornerstone of dynamics). 
furthermore, 
%th namely systems which evolve in time and often display chaotic behavior.    Flows which describe the motion of a particle on a \emph{surface} often arise in  physics (e.g. in celestial mechanics, solid state physics or statistichal mechanics models).  Their study has played a central role in  research in dynamical systems since the birth of the field (with Poincar{\'e}, in the $1890$s). 
many models of systems of physical origin are described by flows on surfaces, starting from celestial mechanics, up to solid state physics or statistical mechanics models.    
The beginning of the study of surface flows can be dated back to Poincar{\'e} \cite{Po:met} at the end of the 19th century,  and coincides with the birth of dynamical systems as a research field.  Poincar{\'e}  was in particular interested in the study of flows on \emph{tori}, or surfaces of \emph{genus one}. Several famous systems in physics %described by by the motion of a point on a \emph{surface} 
lead naturally to the study of flows on surfaces of \emph{higher genus}, which, in this survey, will mean genus $g\geq 2$. Examples include the Ehrenfest model in statistical mechanics (related to a linear flow on a translation surface of genus five),  or the Novikov model in solid state physics, which is described by locally Hamiltonian flows, a class which will be one of the central themes of this survey (see \S~\ref{sec:locHam}).

There is a rich history of results on the topological and qualitative behavior of trajectories (see for example \cite{NZ:flo} and the references therein), as well as on the ergodic theory of certain well studies classes of flows (for example in genus one, in relation with KAM theory, see \S~\ref{sec:g1}, and linear flows on translation surfaces, whose study is intertwined with Teichm{\"u}ller dynamics,  
%and the work of several Fields medallists,
 see \S~\ref{sec:higherg}). Many fundamental problems, though, in particular on the mathematical characterization of chaos (such as dynamical, spectral and rigidity questions) in various natural classes of surface flows,  in particular smooth flows preserving a smooth measure,  were only recently understood  and many others are still  open (see \S~\ref{sec:mixing}).  

One of the reasons for this late development is perhaps that, in order to investigate fine chaotic or rigidity properties of flows in higher genus, one needs to impose quite delicate assumptions on the behavior of orbits on different scales. To capture these multi-scale features, the concept of \emph{renormalization} plays a crucial role (see \S~\ref{sec:renorm}).  
In the case of genus one, the type of assumptions on the flow often take the form of \emph{Diophantine conditions} or, more in general, of \emph{arithmetic conditions}  on the rotation number (see \S~\ref{sec:DChigherg}) and control how well the flow orbits are approximated by \emph{periodic orbits}. The renormalization point of view on these conditions is that they can be described in terms of continued fraction theory and therefore studying the dynamics of the Gauss map, or, equivalently, geometrically, studying the geodesic flow on the modular surface, both of which are classically well understood.  

In higher genus, on the other hand, one had to wait for the development of the rich and fruitful theory of renormalization  in Teichm{\"u}ller dynamics (see \S~\ref{sec:renorm}). This theory provides a renormalization framework (initially developed to study ergodic properties of rational billiard, interval exchange transformations and translation flows), which can be exploited to understand {when} a surface flow is \emph{renormalizable} (see \S~\ref{sec:linhigherg} and \S~\ref{sec:renorm}) and when it preserves a smooth invariant measure; in the latter case, then, it allows  
to impose conditions on a (smooth) surface flows to guarantee the presence of particular chaotic properties (see \S~\ref{sec:locHam}).  
The type and the nature of what we refer to as  \emph{Diophantine-like conditions} in higher genus, which is much more delicate than in genus one and  often involves assumptions on \emph{hyperbolicity} or \emph{uniform hyperbolicity} of the renormalization, will be the leading theme of this survey. These conditions are sometimes also called \emph{arithmetic conditions}, by analogy with the genus one case, even though the relation with  classical arithmetic and Diophantine equations is lost when the genus is greater than one.
 
In what follows, we first start in the next \S~\ref{sec:g1} with the classical case of flows on \emph{genus one surfaces}, recalling some of the classical results on the \emph{linearization problem} and ergodic properties and discussing the related arithmetic conditions. Then, in \S~\ref{sec:higherg}, we  will briefly overview some of the rapid developments in our understanding of ergodic, spectral and disjointness properties of (smooth) area-preserving flows on higher genus surfaces  (see \S~\ref{sec:locHam}), as well as linearization and rigidity problems in higher genus (in \S~\ref{sec:linhigherg}). After having introduced the notion of \emph{renormalization} in this setting  (see \S~\ref{sec:renorm}), we then focus in \S~\ref{sec:DChigherg} on the Diophantine-like conditions behind these results.

\section{Flows on surfaces of genus one and classical arithmetic conditions.} \label{sec:g1}
%The study of flows of genus one surfaces is intimately related to the classical theory of \emph{circle diffeomorphisms}. and to KAM theory (see \S~\ref{sec:circlediffeos}. 
%\subsubsection{Rotations, circle diffeos and linearization}\label{sec:circlediffeos}
A central idea introduced by Poincar{\'e} was that the study of a surface flow can be often \emph{reduced} to the study of a one-dimensional discrete dynamical system, by taking what we nowadays call a \emph{Poincar{\'e}} section and considering  the \emph{Poincar{\'e}} first return map  of the flow to the section (when and where it is defined). If we start from a flow $\varphi_\R:= (\varphi_t)_{t\in\mathbb{R}}$ on a torus, i.e.~on a compact, orientable surface $S$ of genus one and assume that it does not have fixed points,  nor closed orbits (or more generally, nor Reeb components, see \cite{NZ:flo}),  there is a (global) section given by a closed transverse curve and 
 the Poincar{\'e} first return map to it is a  {diffeomorphism} $f:S^1\to S^1$ {of the circle} $S^1\cong \mathbb{R}/\mathbb{Z}$. The simplest example of \emph{circle diffeomorphism} (or \emph{circle diffeo} for short) is a (rigid) \emph{rotation},  i.e.\ the map $R_\alpha(x)=x+\alpha \mod 1$ on $\mathbb{R}/\mathbb{Z}=[0,1]/\sim$.  % The study of flows of genus one surfaces is intimately related to the classical theory of \emph{circle diffeomorphisms}. 
A key concept associated to circle diffeomorphisms is that of \emph{rotation number}: if $\mu$ is an \emph{invariant} probability \emph{measure} for the circle diffeo $f$ (which always exists by Krylov-Bogolyubov theorem), the rotation number $\rho(f)$ of $f$ can be seen as an \emph{average displacement} of points, namely $\rho(f)= \int_{0}^{1} (F(x) -x)\, \mathrm{d}\mu(x)\mod 1$ where $F:\mathbb{R}/\mathbb{Z}\to \R$ is a lift of $f$.  The rotation $R_\alpha$ can be seen as the \emph{linear model}  of a circle diffeo with rotation number $\alpha$. 

The \emph{topological behavior} of trajectories of $(\varphi_t)_{t\in\mathbb{R}}$ can be completely understood and classified exploiting the  \emph{rotation number} (this is essentially the content of \emph{Poincar{\'e} classification theorem}, see \cite{KH:mod} for an expository account): when $\rho(f)\in \mathbb{Q}$,  there exist periodic orbits (which  either \emph{foliate} the surface $S$,  or are attracting or repelling). %, giving rise to  a  \emph{Morse-Smale} type of dynamics). 
% is foliated into periodic orbits for $(\varphi_t)_{t\in\mathbb{R}}$ , or the dynamics is \emph{Morse-Smale}, i.e. there are attracting/repelling periodic orbits). 
On the other hand, when $\rho(f)\notin \mathbb{Q}$, the dynamics of $(\varphi_t)_{t\in\mathbb{R}}$ is either \emph{minimal} on the whole surface (i.e.~all orbits are \emph{dense}), or minimal when restricted to  a \emph{Cantor-like} invariant limit  set (locally product of a Cantor set with $\mathbb{R}$). In the latter case, we speak of \emph{Denjoy-counterexamples}; their existence is ruled out when the diffeo (and the flow) is sufficiently smooth, for example $\mathcal{C}^2$  in view of Denjoy's work \cite{Denjoy} (less regularity, in particular $\mathcal{C}^1$ with bounded variation derivative, suffices,  see e.g.\ \cite{KH:mod} for more details). 
%Among the few \emph{global} results (which do not assume that $T$ is close to $T_0$), we recall that Denjoy showed that as soon as  a circle diffeo $T$ is sufficiently smooth, and  the \emph{rotation number} $\gamma(T)$ of $T$ is \emph{irrational}, $T$ is \emph{minimal} and linearizable (in particular conjugate to the rotation $R_\alpha$ with $\alpha=\gamma(T)$ equal to the rotation number of $T$). 

\subsubsection*{Arithmetic conditions for linearization of circle diffeomorphisms.}
To gain a finer understanding of the dynamics and describe the ergodic behavior of almost-every trajectory with respect to a smooth measure,  one has to address the \emph{linearization problem}, 
a classical question
%problem in dynamics, 
which is at the heart of the {theory of circle diffeomorphisms}. Namely, one wants to understand when a circle diffeomorphism $T$ is \emph{linearizable}, i.e.~conjugate to a rigid rotation $R_\alpha$ (i.e.~when there exists a homeomorphism $h:S^1\to S^1$, called the \emph{conjugacy}, such that $h\circ T=R_\alpha \circ h$) and, if it is linearizable, what is the \emph{regularity} of the conjugacy $h$. 
To address this question, one needs to put further assumptions both on the \emph{regularity} of the diffeo and, in relation to it, the  irrationality of the \emph{rotation number}. 

%\subsubsection*{Classical Diophantine and arithmetic conditions.}
We recall that \emph{arithmetic conditions} are conditions that  prescribe how well (or how \emph{badly}) the irrational rotation number $\alpha\in \mathbb{R}$ is approximated by \emph{rational} numbers and morally control how well the flow orbits are approximated by \emph{periodic orbits}. The best known such condition is perhaps the (classical) \emph{Diophantine condition} (or for short, DC): $\alpha\in \mathbb{R}\backslash\mathbb{Q}$ is said to be \emph{Diophantine} (of exponent $\tau\geq 0$) iff there exists $C>0$ such that
\begin{equation*}
\Big| \alpha -\frac{p}{q}\Big| \geq \frac{C}{q^{2+\tau}}, \qquad \text{for \ all}\ p,q\in\mathbb{Z}, \ q\neq 0.
\end{equation*}
If the above condition holds for $\tau=0$, we say that $\alpha$ is \emph{badly approximable} or \emph{bounded-type}. 
Equivalently, the DC can be rephrased in terms of the continued fraction expansion $[a_0,a_1,\dots, a_n,\dots]$ of $\alpha$: if $q_n$ denotes the \emph{convergents} of $\alpha$, namely the denominators of the partial approximations $p_n/q_n:=[a_0,a_1,\dots, a_n]$, the DC is equivalent to the growth control $a_{n+1}=O(q_n^{\tau})$, while $\alpha$ is of bounded type iff $a_n$ are uniformly bounded. 
% ADD?? See \S~\ref{sec:modular} for an interpretation of these conditions in terms of geodesics on the modular surface. %, via \emph{arithmetic conditions}.

 The \emph{local theory} of linearization of circle diffeos, which treats the case of diffeos $f:S^1\to S^1$ which are $\mathcal{C}^\infty$-\emph{close} (or analytically, or $\mathcal{C}^r$ close) to a circle rotation $R_\alpha$, where $\alpha=\rho(f)$, is a rather classical application of KAM theory. 
%A landmark result is the \emph{local rigidity} theorem of Arnold \cite{Arnold}, who successfully applied KAM theory to show that under a suitable Diophantine-like condition on the rotation number $\rho$, sufficiently small analytic deformations of $x \mapsto x + \rho$, whose rotation number is equal to $\rho$, must be \textit{analytically} conjugate to $x \mapsto x + \rho$. Arnold went on to conjecture that such a statement should hold true without any assumption on the closeness to rotations. 
The prototype result is the  \emph{local rigidity} theorem of Arnold \cite{Ar:sma}, who  showed that if $\alpha$ is Diophantine, circle diffeos which are a 
sufficiently small analytic deformations of $R_\alpha$ and have rotation number equal to $\alpha$, must be \textit{analytically} conjugate to $R_\alpha$. 
% $x \mapsto x + \rho$. 
% $x \mapsto x + \rho$, whose rotation number is equal to $\rho$, must be \textit{analytically} conjugate to $x \mapsto x + \rho$. Arnold went on to conjecture that such a statement should hold true without any assumption on the closeness to rotations. 
%showed that and $f$ smooth, it is smoothly linearizable. 
%Local linearization results require suitable arithmetic conditions related to the regularity of $f$.
Among the few \emph{global results} (which do not assume that $f$ is close to a rotation), we recall the %at Denjoy showed that as soon %as  a circle diffeo $T$ is sufficiently smooth, and  the \emph{rotation number} $\gamma(T)$ of $T$ is \emph{irrational}, $T$ is \emph{minimal} and linearizable (in particular conjugate to the rotation $R_\alpha$ with $\alpha=\gamma(T)$ equal to the rotation number of $T$). 
celebrated theorem by Michael Herman \cite{Herman} and Jean-Christophe Yoccoz \cite{Yo:con}, answering a question by Arnold, that shows that if $f$ is $\mathcal{C}^\infty$ (or analytic) and  its rotation number $\rho(f)$ satisfies the DC, the conjugacy is  
%and   , under a full measure condition  on the rotation number $\gamma(T)$ (which J.C.~Yoccoz showed to coincide with the class of \emph{Diophantine numbers}),% $h$ is
 $\mathcal{C}^\infty$ (resp.~analytic). Furthermore, the DC turns out to be the optimal arithmetic condition for global smooth linearization. Another, more subtle arithmetic  condition, % (which does not depend only on the \emph{asymptotic} growth of entries, but on
 called `\emph{condition} H' in honor of Herman, was introduced by Yoccoz as the optimal condition for global \emph{analytic} linearization of analytic diffeos, see \cite{Yo:CIME}. %We return to this condition and state it in \S~\ref{sec:DChigherg} below.
% Contrary to the DC, which is an \emph{asymptotic condition} on the entries, 
	
%\subsection{Linearization}	

Another famous arithmetic condition is the \emph{Roth-type} condition, which is satisfied by irrationals $\alpha\in\mathbb{R}\backslash\mathbb{Q}$  such that $a_n=O_\epsilon (q_n^\epsilon)$ for all $\epsilon >0$.   A crucial step in the KAM approach developed by Arnold for circle diffeomorphisms is to solve a \textit{linearized} version of the conjugacy equation  $ R_{\alpha}\circ h= h \circ T $, 
%In particular Roth-type numbers are are \emph{Diophantine of all orders}. 
namely the \emph{cohomological equation}:  given a smooth $\phi:I\to \mathbb{R}$, one looks for a smooth solution $\varphi:I \to \mathbb{R}$ to the equation $\varphi\circ R_\alpha- \varphi=\phi$. The Roth-type condition turns out to be the optimal one needed to solve this cohomological equation with optimal loss of differentiability:  for any $r>s+1\geq 1$, one can find a solution $\varphi\in \mathcal{C}^s$ for any $\phi\in \mathcal{C}^{r}$ as long as $\int \phi =0$ (which is a trivial necessary condition), \emph{if and only if} $\alpha$ is Roth-type: this  equivalent characterization provides a remarkable connection between dynamical and arithmetical properties. 

%This equation, known as \textit{cohomological equation}, is easily solved { in the smooth setting}  using Fourier analysis in the case where $R_{\rho} := x \mapsto x + \rho$ is a rotation satisfying a full measure arithmetic condition, under the necessary (and in this setting the only) obstruction that $g$ has to have  zero-mean.
  
%. Such conditions  They usually play a role in \emph{Diophantine-condition}
%, the best know of which is perhaps the , since they relate to how well the \emph{rotation number} is well approximated by rational numbers. 

%\begin{remark}\label{rk:fullmeasure}
We remark that the Diophantine condition, the H condition and the Roth type condition can all be proved to have \emph{full measure}, namely they hold for a set of $\alpha\in [0,1]$ of Lebesgue measure one (the set of badly approximable $\alpha\in[0,1]$, on the other hand, has Lebesgue measure \emph{zero}, although  full Hausdorff dimension). While full measure of the Diophantine and Roth conditions can be proved in an  elementary way, it is an instructive exercise to derive it
%These results can be proved %as (simple) ergodic theory exercises
% using 
from the properties of the Gauss map $G:[0,1]\to [0,1]$ and of the Gauss  invariant measure  ${\mathrm{d}x}/{{\log 2}(1+x)}$, since this point of view can be applied to show full measure of  other arithmetic conditions as well and it can be furthermore  generalized to higher genus (see \S~\ref{sec:renorm} and \S~\ref{sec:DChigherg}).  
% $\frac{\mathrm{d}x}{\sqrt{2}(1+x)}$ 
%preserved by $G$, whose symbolic coding produces the continued fraction entries of $\alpha$. 
%\end{remark}
%\noindent 

In view of this remark, we conclude this section  with a reinterpretation of Herman's linearization theorem in the language of \emph{foliations}  into flow trajectories. 
%Throughout this survey, the word \emph{typical} will be always used in the measure theoretical sense.
\begin{theorem}[reformulation of Herman's global theorem \cite{Herman}]\label{HYthm}
For a full measure set of real numbers $\alpha$,  a foliation  on a genus one surface which is topologically  conjugate to the foliation given by a linear flow with rotation number $\alpha$, is also $\mathcal{C}^\infty$ conjugate to it.
\end{theorem}

\subsubsection*{Ergodic properties in genus one and exceptional behavior.}
From the existence (and abundance) of smooth (or at least  continuously \emph{differentiable}, i.e.~$\mathcal{C}^1$) linearizations, one can infer many of the smooth  measure theoretical  ergodic properties of flows in genus one.    
%\begin{theorem}
%For a \emph{typical} rotation number,  a smooth flow $\phi_\mathbb{R}$ on a surface of genus one is uniquely ergodic with respect to a measure with an absolutely continuous invariant measure 
%\end{theorem}
In particular, one sees that, for a full measure set of rotation numbers, flows in genus one are \emph{ergodic} (since irrational rotations are)  with respect to a \emph{smooth} invariant measure of full support  (the $\mathcal{C}^1$-regularity of the conjugacy allows indeed to \emph{transport} the Lebesgue invariant measure to obtain the invariant measure for the diffeo, which in turns give a \emph{transverse measure} for the flow). Furthermore, they are \emph{uniquely ergodic}  (in view of Kronecker-Weyl theorem for rotations, e.g.~\cite{CFS:erg2}), i.e. this natural invariant measure is the \emph{unique} invariant measure (up to scaling). %This implies in particular that one can control the asymptotic behavior of \emph{all} trajectories.  
%To infer further ergodic properties of the flow from its Poincar{\'e} map (in particular properties of spectral nature, see \cite{}, which depend not only on the trajectories of the flow, but also on its parametrization), one can exploit the \emph{special flow} representation, i.e.~take into consideration also the \emph{return time} to the section. We explain this contruction and return to this point in \S~\ref{sec:specialflows}. %One cna show in particular that no flow in genus one is \emph{mixing}.
%: a flow $\varphi_\R$ with Poincar{\'e} map $T:I\to I$ preserving a measure $\nu$ is (measure theoretically) isomorphic to a flow on ADD

%Smooth flows 
%for typical rotation numbers, (i.e.~for a full measure set of $\alpha$), a smooth flow with rotation number $\alpha$ is \emph{ergodic
%since sufficiently regular flows on tori are  conjugate to the linear counterparts (see the previous section), they display quite a rigid dynamics and faible ergodic properties. 
We remark that  \emph{exceptional} ergodic behaviors in genus one ({smooth}) surface flows, can be constructed  for  flows  whose rotation numbers are irrational but not Diophantine, i.e.~so called \emph{Liouvillean} (rotation) numbers. When $\alpha$ is Liouville, exploiting the abundance of good rational approximations $p_n/q_n$ to $\alpha$, for example using the method of \emph{periodic approximations} pioneered by Anosov  and Katok 
and  later revived  by Fayad, Katok  et al (see  \cite{KH:mod} or  the survey \cite{FK:Lio}), one can construct many examples  with \emph{pathological} behavior, for example flows with a \emph{singular} invariant measures and \emph{time-reparametrizations} (also called  \emph{time-changes}) which are weakly mixing or which have mixed spectrum (see %\cite{KH:mod, FKW:mix} as well as 
\cite{FK:Lio} and the references therein). 

%\subsubsection*{Flows with singularities in genus one}
Finally, before moving to higher genus, we remark that another possible way to introduce interesting dynamical features for \emph{typical} rotation numbers, is to consider flows on tori \emph{with singularities}. The simplest type of singularity is a \emph{stopping point}. Already such a simple perturbation, which is only a time-reparametrization  of the flow, can lead to  flows which are typically  \emph{mixing} (see \cite{Ko:mix}) and even to flows with \emph{Lebesgue spectrum} (see \cite{FFK:Leb}).  %simplest example are flows with one simple saddles and one center .  
Smooth measure preserving flows on a torus with one center and one simple saddle (see Figure~\ref{Arnoldtorus}) were first studied by Arnold in  \cite{Ar:top} and constitute one of the  most studied examples in the class of flows known as \emph{locally Hamiltonian}: 
%(which will be defined in the next \S~\ref{sec:locHam}).  
we return to them and to their typical ergodic properties in \S~\ref{sec:mixing}.

%The simplest examples of locally Hamiltonian flows  with singularities on a torus, i.e.~flows with  (see Figure~\ref{Arnoldtorus}), were studied by V.~Arnold in \cite{Arn} and are nowadays often called \emph{Arnold flows}\footnote{More precisely, referring to the decomposition described in \S~\ref{sec:generic}, we call \emph{Arnold flow} the restriction to a minimal component obtained by removing the center and the disk filled by periodic orbits around it (called \emph{island}), which, as Arnold shows in  \cite{Arn}, is always bounded by a saddle loop.\label{Arnoldflowfootnote}}.
% is to consider flows on tori which have both a center and a simple saddle. These flows, known as \emph{Arnold flows}, were studied by Arnold %in \cite{Arn:top}, who showed that conjectured that typically they are mixing, as later proved by Khanin and Sinai \cite{SK:mix} (see also %\cite{Ko:07} for the optimal arithmetic condition).
%These flows admit a representation as \emph{special flows} 

%Another possibility way to introduce interesting dynamical features for \emph{typical} rotation numbers, is to consider flows on tori \emph{with singularities}. We return to this below, see \S~\ref{sec:singularities}).

\section{Dynamics of flows on surfaces of  higher genus.}\label{sec:higherg}
Let us now consider the \emph{higher genus} case, namely consider now a (smooth)  
flow $\varphi_\R:= (\varphi_t)_{t\in\mathbb{R}}$  on a  compact, connected orientable (closed) surface $S$ of genus $g\geq 2$. 
%Let $M$ be a compact, connected, orientable  closed (smooth) surface and let $g$ denote its genus. 
Notice that in this case, by Euler characteristic restrictions, the flow \emph{always} has \emph{fixed points} (see Figure~\ref{locHamsing} for some examples). %Remark that when $g\geq 2$,  $Fix(\varphi_\R)$ is always not empty; 
We  require that  singularities be \emph{isolated} (so that in particular, by compactness, the set 
$Fix(\varphi_\R)$ of fixed points is \emph{finite}).

%\subsection{Topological structure and basic constructions}

\subsubsection*{Topological dynamics and quasi-minimal sets}
The \emph{topological classification} of the possible behavior of trajectories of a flow on a surface (and more in general of surface \emph{foliations} which are not necessarily orientable) has been a topic of research in the $20$th century (starting from the $1930s$-$40s$, up to the $1970s$).  In particular, through the works of Maier, %\cite{Ma:tra}, 
Levitt,  %\cite{Le:feu, Levitt2}, 
Gutierrez, Gardiner et al (see \cite{NZ:flo} for references) one could obtain results on what are possible \emph{orbit closures}, as well as  a classification of \emph{quasi-minimal sets}, which can be defined as possible $\omega$-\emph{limit sets} of \emph{non-trivial} recurrent trajectories, i.e.~set of accumulation points of trajectories different from a fixed point or a closed, periodic orbit. Quasi-minimal sets can be the whole surface, subsurfaces with boundary, or a Cantor-like invariant sets. Moreover, one can prove  \emph{decomposition theorems} that shows that one can \emph{cut}  the surface $S$ into subsurfaces 
%which play the role of \emph{basin of attraction} for the
each of which contains at most one quasi-minimal set (see in particular the work by Levitt \cite{Le:dec}). We do not enter here in the details of these topological results, but refer the interested reader for example to the monograph \cite{NZ:flo} and the references therein.

%\subsection{Interval exchanges and generalized IETs as Poincar{\'e} sections}
\subsubsection*{Interval exchanges and generalized IETs as Poincar{\'e} sections} 
As in the case of genus one, an essential tool to study a higher genus flow is to consider a (local) \emph{transversal}  $I\subset S$ to the flow %(which exists as soon as there is a non-trivial recurrent trajectory, i.e.~a recurrent trajectory which is not a fixed point nor a closed orbit, see e.g.~\cite{NZ:flo}) 
and  the \emph{Poincar{\'e} first return map} $T$ of the flow on $I$ (when it is defined, for example almost everywhere when the flow preserves a  finite measure with full support; see more generally  \cite{NZ:flo}). 
%, the existence of a global transversal is guaranteed by Poincar{\'e} Bendixon theory.  
Such first return maps $T:I\to I$  are one-to-one \emph{piecewise diffeomorphisms} known as \emph{generalized interval exchange transformations}:   a map $T:I\to I$  is a generalized interval exchange transformations or, for short, a GIET, if one can partition $I$ into  intervals $I_1, \dots, I_d$ (finitely many since we are assuming that $\varphi_\R$ has finitely many fixed points) so that the restriction $T_i$ of $T$ to $I_i$, for each $1\leq i\leq d$, is a diffeomorphism onto its image which extends to a diffeo of the closure $\overline{I}_i$ (see e.g.~\cite{MMY:lin}). We say in this case that $T$ is a $d$-GIET. We say furthermore that $T$ is \emph{of class} $\mathcal{C}^r$ if the restriction of $T$ to each $I_i$ extends to a $\mathcal{C}^r$-diffeomorphism onto the closed interval $\overline{I}_i$. 
The adjective \emph{generalized} is used to distinguish them from the more commonly studied (standard) interval exchange transformations (or simply IETs), which are one-to-one piecewise \emph{isometries}, 
% so that their graph if \emph{piecewise linear} with derivative one 
namely GIETs such that the derivative $T_i'$ of each branch is  constant and equal to one.  %(so that $T_0$ is a piecewise isometry, 
%(see Figure \ref{IET}, right).  

%\begin{figure}[h]
%\subfloat[GIET \label{GIET}]{\includegraphics[width=3.6cm]{Gietgraph}}
%\qquad
%\subfloat[IET\label{IET}]{\includegraphics[width=4.1cm]{IETplot}}
%\caption{A generalized  IET (GIET) and a (standard) IET with $d=4$.\label{IET}. \label{locHamflows}}
%\end{figure}
%Notice that rigid rotations are IETs with $d=2$ branches. 
Standard IETs are a generalization of circle rotations (since a IET is a rotation when $d=2$) and play an analogous role in higher genus, providing the natural \emph{linear model} of a GIET (see \S~\ref{sec:linhigherg}). 
%namely GIETs such that the derivative $T_i'$ of each branch is  constant and equal to one (so that $T_0$ is a piecewise isometry, see Figure \ref{IET}, right). Notice that rigid rotations are IETs with $d=2$ branches.
Furthermore, as rotations are Poincar{\'e} maps of \emph{linear flows} on the torus $\mathbb{R}^2/\mathbb{Z}^2$ (i.e.~flows which arise as solutions of $(\dot{x_1}, \dot{x_2})=(\theta_1,\theta_2)$, which move points with unit speed along lines of slope $\theta_2/\theta_1$), 
IETs arise naturally as Poincar{\'e} maps of \emph{linear flows} on \emph{translation surfaces} (see the ICM proceeding \cite{CW:ICM} for an introduction to the latter).

\subsection{Locally Hamiltonian flows}\label{sec:locHam}
%\subsubsection*{Locally Hamiltonian flows and their special flow representation}
We will be mostly concerned with flows which preserve a (probability) measure $\mu$ of \emph{full support}, for example an area-form, since this is a natural setup for \emph{ergodic theory}. Given a surface $S$ with a fixed smooth area form  $\omega$, a \emph{smooth area preserving flow}     $\varphi_\mathbb{R}=(\varphi_t)_{t\in\mathbb{R}}$ on $S$ is  a smooth flow on $S$ which preserves the measure $\mu$ given  integrating a smooth density with respect to $\omega$. The interest in the study of these flows  and, in particular, in their ergodic and mixing properties, was revived by Novikov \cite{No:the} in the $1990$s, in connection with problems arising in solid-state
physics as well as in pseudo-periodic topology (see e.g.\ the survey \cite{Zo:how} by A.~Zorich).  Smooth area-preserving flows are  also called \emph{locally Hamiltonian flows} or \emph{multi-valued Hamiltonian flows} in the literature, in view of their interpretation as flows locally given by Hamiltonian equations:  
one can find local coordinates $(x_1,x_2)$ on $S$ in which $\varphi_\mathbb{R}$ is given by 
the solution to the  equations: 
 $$\begin{cases}\dot{x_1}&={\partial H}/{\partial x_2},\\ \dot{x_2}& =-{\partial H}/{\partial x_1},\end{cases}$$
 where $H$ is a   real-valued (\emph{local}) Hamiltonian. For simplicity we will assume here that $H$ is infinitely differentiable, even though for several results $\mathcal{C}^3$ (or also $\mathcal{C}^{2+\epsilon}$ for every $\epsilon>0$)  suffices.
% A \emph{global}  Hamiltonian $ H$ cannot be in general be defined (see \cite{NZ:flo}, Section 1.3.4), but one can think of  $\varphi_\mathbb{R}$ as globally given by a \emph{multi-valued} Hamiltonian function. 
%\footnote{Novikov \cite{No} and his school in the $1990s$  advocated the study of locally Hamiltonian flows as model to describe the motion of an electron in a metal under a magnetic field in the semi-classical approximation (the surface appears here as Fermi energy level surface).}
It turns out that such smooth area preserving flows on $S$  are in one-to-one correspondence  with smooth \emph{closed} real-valued differential $1$-forms: given such a $1$-form $\eta$, we can associate to it the integral flow $\varphi^\eta_\R$ of  the vector field  $X$ such that $\eta = i_X \omega$, where $i_X$ denotes the contraction operator.
%, i.e. $i_X \omega =\omega( \eta, \cdot )$ and consider the flow $\varphi_\mathbb{R}$ on $S$ given by $X$. 
Since $\eta$ is closed, $\varphi^\eta_\R$ is  area-preserving; conversely, every smooth area-preserving flow can be obtained in this way.

\begin{figure}[h]
\subfloat[\label{Arnoldtorus}]{\includegraphics[width=3.96cm]{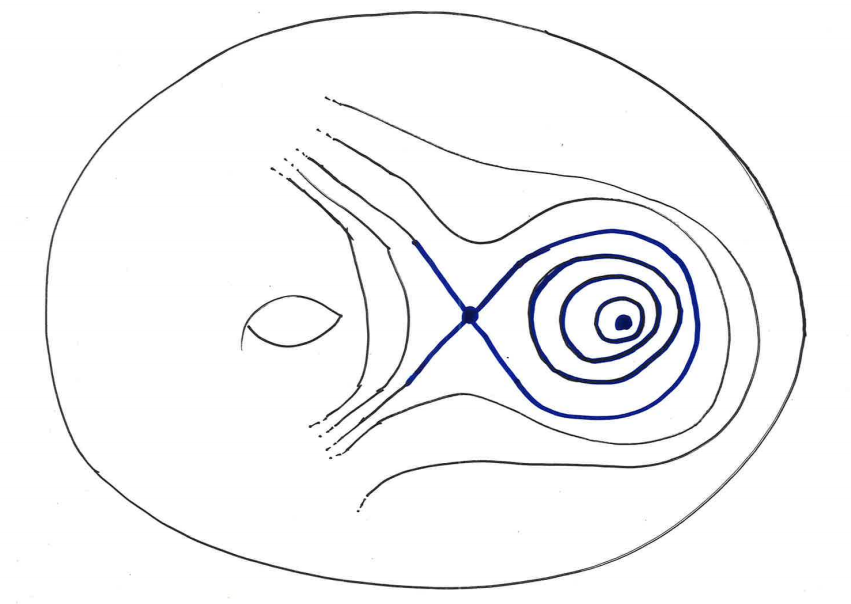}}\hspace{3mm}
\subfloat[\label{g3}]{\includegraphics[width=79mm]{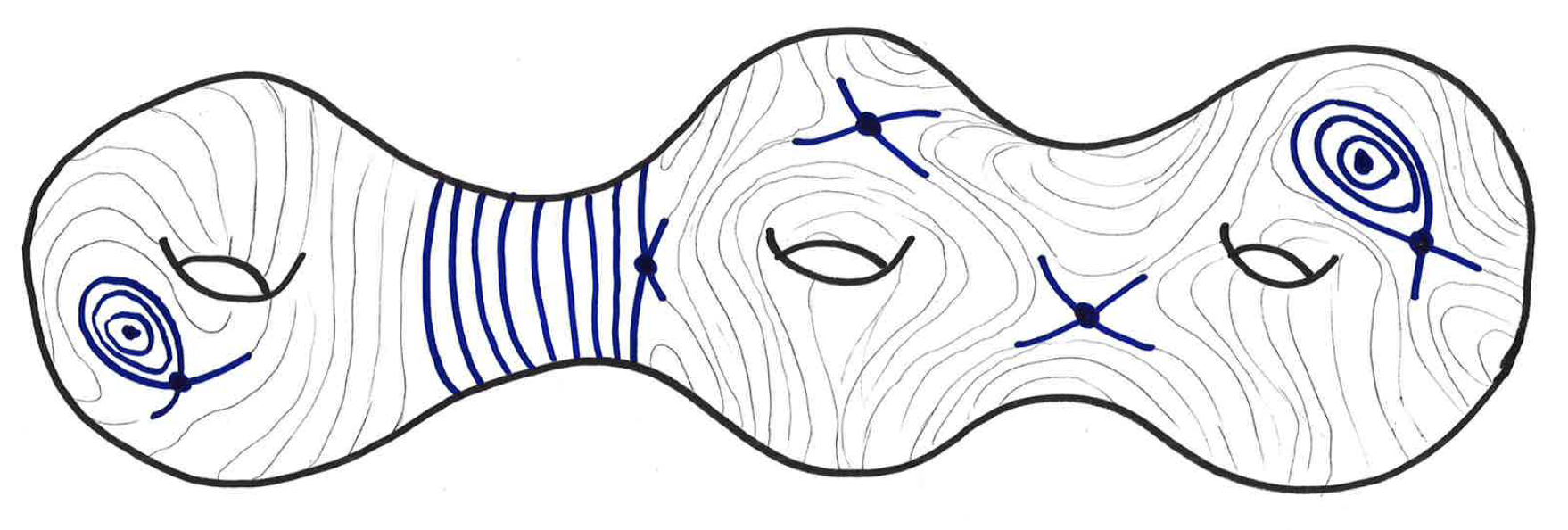}}
\caption{Pictorial representation  of  locally Hamiltonian flows on a surfaces: in \textbf{(a)}  an Arnold flow ($g=1$) and in \textbf{(b)}  a flow in $g=3$ with two minimal components and $3$ periodic components. \label{locHamflows}}
\end{figure}

%\subsubsection*{Special flow representation}
%ADD definition of locally Hamiltonian flows 
%If the measure is smooth, 
%%\subsubsection*{Special flow representation}
%A standard construction in ergodic theory allow to recover a flow  $\phi_\mathbb{R}$ (up to measure theoretical isomorphism) from the its %Poincar{\'e} map \emph{and} the knowledge of the \emph{return time function}, namely a function $r:I\to \mathbb{R}^+$ which gives the \emph{first return time} $r(x)$ of the flow-orbit of $x\in I$ under $\phi_\mathbb{R}$
%Linear flows on translation surfaces, for example, can be represented as special flows over IETs, under a roof function $r$ which is \emph{piecewise constant}, since all points in the same continuity interval for $T$ return at the same time. On the other hand, locally Hamiltonian flows If the measure is \emph{smooth} (i.e.~it is obtained by integrating a smooth  % $f(x,y)\omega$ where $f$ is smooth and $\omega $ is a 
%non-degenerate two form on $S$), and hence induces a smooth (in particular non atomic) \emph{transverse measure} on the foliation into flow trajectories, one can choose suitable coordinates such that the Poincar{\'e} map is a \emph{standard} IET.  

\subsubsection*{Topology and measure class}
\begin{figure}[h]
\subfloat[%center 
\label{center}]{\includegraphics[width=2.6cm]{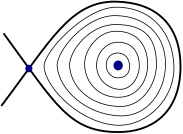}}\hspace{8mm}
\subfloat[%simple saddle 
\label{simplesaddle}]{\includegraphics[width=2cm]{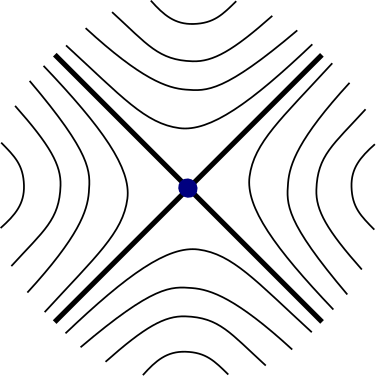}}\hspace{8mm}
\subfloat[%multisaddle 
\label{multisaddle}]{\includegraphics[width=1.8cm]{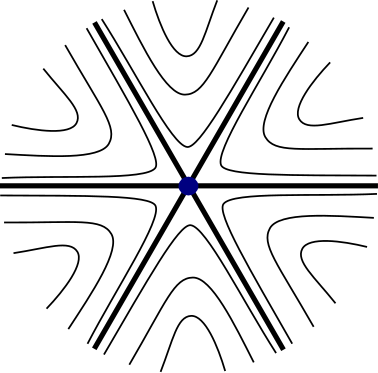}}\hspace{8mm}
%\qquad \quad
\subfloat[%multisaddle 
\label{shearing}]{\includegraphics[width=2cm]{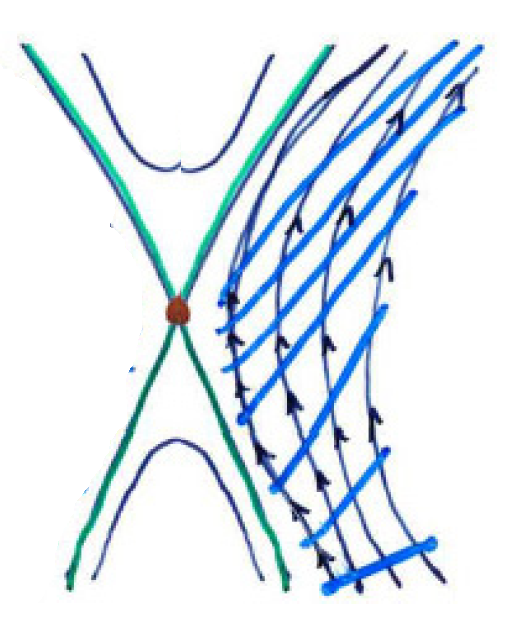}}
\caption{Type of singularities of a  locally Hamiltonian flow:  a center in \textbf{(a)}, a simple saddle in \textbf{(b)} and a multisaddle in \textbf{(c)}. Decelerations and shearing near an Hamiltonian  saddle in \textbf{(d)}. \label{locHamsing}}
\end{figure}

Let $Fix(\varphi_\R)$ denote the set of \emph{fixed points} (also called \emph{singularities}) of the flow $\varphi_\R$. Remark that when $g\geq 2$,  $Fix(\varphi_\R)$ is always not empty; we  require  the singularities to be \emph{isolated} (so that in particular, by compactness, 
$Fix(\varphi_\R)$ is a \emph{finite set}) and  denote by $\mathcal{F}$ the set of smooth closed $1$-forms on $S$  with isolated zeros. 
On $\mathcal{F}$ (which we can think of as the space of locally Hamiltonian flows) one can define a \emph{topology} as well as a  measure class. The  \emph{topology} %on $\mathcal{F}$ 
is obtained by considering perturbations of closed smooth $1$-forms by (small) closed  smooth $1$-forms.  
%\footnote{Let $\eta$, $\eta'$ be two smooth closed $1$-forms.  We say
%that $\eta'$ is an $\evarpsilon$-perturbation of $\eta$ if for any $x\in S$ there exists coordinates on a simply connected neighbourhood  $U$ of $x$, such that  $\eta \Vert _U=dH $ and $(\eta'-\eta)\Vert_U= dh$ where $\Vert h\Vert_{\infty}<\epsilon \Vert H\Vert_{\infty}$ (here $\Vert \cdot \Vert_{\infty}$ denotes the $\mathcal{C}^{\infty}$ norm).}. We say that a condition is \emph{generic} (in the sense of Baire) if it holds for
%flows described by an open and dense set of forms with respect to this topology.
We will often restrict our attention to the subset $\mathcal{M} \subset \mathcal{F}$ of \emph{Morse} closed $1$-forms, % (adopting the notation introduced by Ravotti \cite{Ra:mix}), 
(i.e.\ forms which are locally the differential of a \emph{Morse function}), which is \emph{open and dense}  in $\mathcal{F}$ with respect to this topology (see e.g.\ \cite{Ra:mix}).  Locally Hamiltonian flows corresponding to forms in $\mathcal{M}$ have only \emph{non-degenerate fixed points}, i.e.\ \emph{centers}  and \emph{simple saddles} (as in Figures~\ref{center} and ~\ref{simplesaddle}), as opposed to degenerate \emph{multi-saddles} (as in Figure \ref{multisaddle}).

%which only have \emph{non-degenerate singularities}, so in particular only simple saddles and centers,   
 %(i.e.~ a function that has \emph{non-degenerate} zeros, so that the Hessian at every fixed point is non degenerate). The set $\mathcal{A}$ of Morse  $1$-forms 
%is \emph{open and dense}  in $\mathcal{F}$ (with respect to the topology defined above). 

 A \emph{measure-theoretical notion of  typical} can be defined  using the \emph{Katok fundamental class} (introduced by Katok in \cite{Ka:inv}, see also \cite{NZ:flo}), i.e.\ the cohomology class of the 1-form $\eta$  which defines the flow. More explicitely, if we fix a base 
 $\gamma_1, \dots, \gamma_n$ of the relative homology $H_1(S, Fix(\varphi_\R), \mathbb{R})$ (where $n=2g+k-1$ if $k$ denotes  the cardinality of  $Fix(\varphi_\R)$) and consider the period map $Per $ given by  $Per(\eta) = (\int_{\gamma_1} \eta, \dots, \int_{\gamma_n} \eta) \in \mathbb{R}^{n}$,  we say that a property holds for a \emph{typical} locally Hamiltonian flow in $\mathcal{F}$, if it holds for all $\eta$ such that $Per(\eta)$ belongs to a full measure set with respect to the Lebesgue measure on $\mathbb{R}^n$.

\subsubsection*{Minimal components and ergodicity.}
To describe (typical) chaotic behavior in locally Hamiltonian flows, it is crucial to distinguish between two  open sets (complementary, up to measure zero, see \cite{Ul:slo} or \cite{Ra:mix} for more details): in the first open set, which we will denote by $\mathcal{U}_{min}$, the typical flow is \emph{minimal} (the term \emph{quasi-minimal} is also used in the literature), in the sense that the orbits of all points which are not fixed points are \emph{dense} in $S$;  flows in $\mathcal{U}_{min}$ have only saddles, since the presence of centers prevents minimality. 
On the other open set, that we call $\mathcal{U}_{\neg min}$, the flow is not minimal (there are saddle loops homologous to zero which disconnect the surface), but one can decompose the surface into a finite number of subsurfaces with boundary $S_i$, $i=1,\dots ,N$ such that for each $i$ either $S_i$ is a \emph{periodic component}, i.e.~the interior of $S_i$ if foliated into closed orbits of $\varphi_\mathbb{R}$ (in  Figure~\ref{locHamflows} (b) one can see three periodic components, namely two disks and one cylinder), or $S_i$ is such that  the restriction of  $\varphi_\mathbb{R}$ to $S_i$ is minimal in the sense above, as pictured in the remaining two subsurfaces in Figure~\ref{locHamflows} (b). These are called \emph{minimal components} and there are at most $g$ of them (where $g$ is the genus of $S$), see \S~\ref{sec:mixing}. 
% (A characterizing property of  $\mathcal{U}_{\neg min}$  is the abscence of \emph{saddle connections} (namely trajectories which join a saddle to a saddle) which are \emph{homologous to zero} (i.e.~split the surface in two disjoint components. While other saddle connections can be removed by perturbation, those homologous to zero form an open set. The set $\mathcal{U}_{\neg min}$  is an open and dense subset of locally Hamiltonian flows \emph{with} saddle loops homologous to zero)

Notice that minimality and ergodicity of a (minimal component of a) locally Hamiltonian flow  are equivalent to minimality or respectively ergodicity of an (and hence any) interval exchange transformation which appears as the Poincar{\'e} map. Classical results proved in the $1980s$  guarantee that almost every IET (with respect to the Lebesgue measure on the interval lengths, assuming that the permutation is irreducible) is minimal (as showed by Keane \cite{Ke:int}, see also \cite{Ka:inv}) and (uniquely) ergodic (as proved in the  works by Masur \cite{Ma:erg} and Veech \cite{Ve:gau}, considered early milestones of the successful application of Teichm{\"u}ller dynamics to the study of IETs and translation surfaces, see the ICM proceeding \cite{CW:ICM} or the survey \cite{Zo:fla}).   
It then follows from definition of Katok measure class that  a typical local Hamiltonian flow in $\mathcal{U}_{min}$ is minimal and  ergodic and, given a  typical local Hamiltonian flow in $\mathcal{U}_{\neg min}$, its restriction on each minimal component is ergodic.

\subsubsection*{Classification of mixing properties}\label{sec:mixing}
%In the classification of mixing properties in particular,  it is crucial to distinguish between two  open sets (complementary, up to measure zero, see \S~\ref{sec:ergodicity} for more details): in the first open set, which we will denote by $\mathcal{U}_{min}$, the typical flow is \emph{minimal}, in the sense that the orbits of all points which are not fixed points are \emph{dense} in $M$. On the other open set, that we call $\mathcal{U}_{\neg min}$, 
Finer chaotic features of locally Hamiltonian flows, in particular mixing and spectral properties,  change according to the type of singularities  and depend crucially on  the locally Hamiltonian parametrization of saddle points.  For a (non-generic) locally Hamiltonian flow with at least one \emph{degenerate saddle} (e.g.~as in Fig.~\ref{multisaddle}), \emph{mixing} (see for the definition \eqref{eq:multimixing} for $n=2$) was proved in the $1970s$ by Kochergin  \cite{Ko:mix}. When on the other hand $\eta\in \mathcal{M}$ is a Morse-one form, so that all saddles are \emph{simple}, one has a dichotomy: inside the open set $\mathcal{U}_{min}$ in which the typical flow is minimal, almost every locally Hamiltonian flow is \emph{weakly mixing}, but it is \emph{not mixing}  in view of work by the author \cite{Ul:wea, Ul:abs}; on the other hand, for a full measure set of flows in  $ \mathcal{U}_{\neg min}$,  the restriction to each of minimal components is mixing (as proved by Ravotti \cite{Ra:mix} extending previous work by the  author \cite{Ul:mix}). 

The question of mixing in higher genus was raised  by V.\ Arnold in the $1990$s, when he conjectured (see \cite{Ar:top}) that  the restriction of a typical smooth flow on a torus with one center and one simple saddle to its minimal component (namely for what we nowadays call an \emph{Arnold flow}) was indeed mixing. His conjecture was proved shortly after  by Khanin and Sinai in \cite{SK:mix}, who showed mixing under the assumption that the rotation number $\alpha$ is such that the entries $a_n$ of the continued fraction expansion of $\alpha$ do not grow too fast, namely there exists a power $1< \tau< 2$ and $C>0$ such that $|a_n|\leq C n^{\tau}$. One can show (for example exploiting the Gauss map $G$ and the finiteness of 
$\int_0^1 a_0(x) \, \mathrm{d}\ \mu_G (x)$,  where $\mu_G$ is the Gauss measure,  
via a standard Borel-Cantelli argument) that this arithmetic condition holds for a  full measure set of $\alpha$. The condition was later improved by Kocerghin, see \cite{Ko:som}. Also in the case of absence of mixing, a prototype result for flows over a full measure set of rotation numbers  was proved by Kochergin \cite{Ko:abs} already in the $1970s$ (and much more recently extended in \cite{Ko:07} to all irrational rotation numbers), much earlier than results in higher genus \cite{Sch:abs, Ul:the, Ul:abs}. 

In higher genus, the above mentioned results on mixing/absence of mixing require the introduction of Diophantine-like conditions, which describe the full measure set of locally Hamiltonian flows for which the results hold. In \cite{Ul:mix}, for example, we introduced a condition on a IET (see \S~\ref{sec:DChigherg} for more details) called \emph{Mixing Diophantine Condition} (or MDC for short). Let us say that the restriction of a  locally Hamiltonian flow $\varphi_\R$ to one of its minimal components $S_i$ satisfies the MDC if one can find a  section $I\subset S_i$ (in \emph{good position} in the sense of \cite{MMY:lin}) such that the IET which arise as Poincar{\'e} map of $\varphi_\R$ to $I$ satisfies the MDC. 
One can then prove:
\begin{theorem}[Ulcigrai \cite{Ul:mix}, Ravotti \cite{Ra:mix}]\label{thm:mixing}
Let $\varphi_\R$ be a flow in $ \mathcal{U}_{\neg min}$ and let $S_i$ be a minimal component. If the restriction of   $\varphi_\R$ to $S_i$ satisfies the Mixing Diophantine Condition, then  $\varphi_\R$ restricted to $S_i$ is mixing.
\end{theorem}
\noindent We then show in \cite{Ul:mix} (exploiting results from \cite{AGY:exp}, see \S~\ref{sec:DChigherg}) that the MDC is satisfied by a full measure set of IETs. 
Similarly, to prove that a typical flow in $\mathcal{U}_{min}$ is \emph{not} mixing, a Diophantine-like condition, which is later proved to be of full measure, is introduced in  \cite{Ul:abs}. Special cases of the absence of mixing result for surfaces with $g=2$ and  two isometric saddles were proved in \cite{Ul:the} and by Scheglov in \cite{Sch:abs}.
% by the second author  (see also \cite{Ko:abs, Ko:abs2} and \cite{Sch:abs} for previous special cases of this result). 
We remark that in $\mathcal{U}_{min}$ there exist nevertheless  exceptional mixing flows, as shown by the work by \cite{ChW:mix}, which produces sporadic  examples in $g=5$.

\subsubsection*{Parabolic dynamics and slow chaos}
Smooth, area-preserving flows on surfaces also provide  one of the fundamental classes of \emph{parabolic}, or \emph{slowly chaotic}, dynamical systems (see e.g.\ the survey \cite{Ul:slo}).  
In  systems which display \emph{sensitive dependence} on initial conditions (the so-called   \emph{butterfly effect}), one can find many nearby initial conditions whose trajectories \emph{diverge} with time. Contrary to hyperbolic systems, where this divergence happens (infinitesimally) at \emph{exponential speed}, in parabolic systems the divergence speed is \emph{slow}, namely sub-exponential and in all known examples \emph{polynomial} or  sub-polynomial. Slow divergence in locally Hamiltonian flows  is created by Hamiltonian saddles, which create different deceleration rates of nearby trajectories and produce a form of (local) \emph{shearing}, by \emph{tilting} in the flow direction the image under the flow of arcs initially transverse to the dynamics, as illustrated in Figure~\ref{shearing}. Shearing happens not only locally, near a saddle, but  globally for typical flows in $ \mathcal{U}_{\neg min}$, which (in view of the presence of saddle loops) display a global  \emph{asymmetry} in the prevalent direction of shearing. It is this geometric mechanism which is behind the proof of mixing (in this setting, but also for many other classes of parabolic flows, see the survey \cite{Ul:she} and the references therein). Under the assumption that the restriction of $\varphi_\R$ to  a minimal component $S_i$ satisfies the Mixing Diophantine Condition, one can produce quantitative estimates on shearing of transverse arcs and, as shown by Ravotti in \cite{Ra:mix}, prove quantitative mixing estimates, which show that mixing happens (at least) at \emph{sub-polynomial speed}, i.e.\ for any two smooth observables $f,g:S_i\to \R$  supported outside the saddles in $Fix (\varphi_\R)\cap S_i$, 
$$
\left |\int_{S_i} f(\varphi_t(x)) g(x) \mathrm{d} \mu - \int_{S_i} f \mathrm{d} \mu \int_{S_i} f \mathrm{d} \mu \right| \leq \frac{C_{f,g}}{(\log t)^\gamma}, \qquad t\geq 0.
$$
This is expected to be also the optimal nature of the estimates, namely the decay is \emph{not} expected to be polynomial or faster in this setting, but no lower bounds on the decay of correlations are currently available. 
%A key component in our project will consist of exploiting the description of quantitative forms of parabolic shearing  to study fine chaotic and spectral properties. 

%Ravotti also shows in \cite{Ra:mix} \emph{subpolynomial} bounds for the speed of mixing.

%Further recent work (see \cite{KKK}) also shows that locally Hamiltonian flows in $\mathcal{U}_{\neg min}$
%display a \emph{quantitative shearing} property inspired by thewhich plays a crucial role in the theory of unipotent flows (ore more precisely a variation introduced in ~\cite{FK} to deal with the presence of singularities).
% From this property, one can deduce  that the restriction of a typical locally Hamiltonian flow $\varphi_\R$ in $\mathcal{U}_{\neg min}$  on its minimal components is not only mixing, but \emph{mixing of all orders}, see \cite{KKU}. 
 %\subsubsection*{Ratner-like properties}

%\begin{figure}[h]
%\includegraphics[width=9cm]{ShearingArc}
%\caption{Transverse segment shearing in the direction of the flow. }
%\end{figure}

\subsubsection*{Ratner's forms of shearing.}
Striking consequences of shearing (such as measure and joining rigidity) were proved for another famous class of parabolic flows, namely horocycle flows on hyperbolic surfaces and their time-changes, by exploiting a \emph{quantitative shearing} property introduced by Marina Ratner and nowadays known as  \emph{Ratner property} (or RP). In view of its importance in the study of horocycle flows and more generally unipotent flows in homogeneous dynamics, it is natural to ask whether this property can be proved and exploited in other parabolic (non homogeneous) settings. For locally Hamiltonian flows, which are natural candidates, the original Ratner property is believed to fail due to the presence of singularities (see \cite{FK:mul}). Nevertheless, a variant of the RP which has the same dynamical consequences, called \emph{Switchable Ratner Property} (or SRP for short),  was introduced by B.\ Fayad and A.~Kanigowski in \cite{FK:mul} and showed to hold for typical Arnold flows (as well as some flows in genus one with one degenerate singularity). As an abstract consequence of the SRP property, one can conclude that typical Arnold flows are not only mixing, but \emph{mixing of all orders}, namely for any $n\geq 2$ and any $n$-tuple $A_0$,$\dots$, $A_{n-1}$ of measurables sets,
\begin{equation}\label{eq:multimixing}
\mu \left(A_0 \cap \varphi_{{t_1}}({A_1})\cap %\varphi_{{t_1+t_2}}({A_2}) \cap 
\cdots \cap \varphi_{{t_1+\dots +t_{n-1}}}({A_{n-1}}) \right)
 \xrightarrow{t_1,t_2\dots, t_{n-1} \to \infty}   \mu( { A_0}) \cdots \mu ({A_{n-1}} ).\
\end{equation}
Notice that this definition reduces to the classical definition of mixing in the special case $n=2$;  whether mixing implies mixing of all orders in general is still an open problem, known as \emph{Rohlin conjecture}.% which this result confirms in this setting. 
 
To prove the SRP property, one needs to assume that the rotation number $\alpha=[a_0, a_1,\dots, a_n, \dots ]$ satisfies an ad-hoc arithmetic condition, namely, if $q_n$ are the denominators  of $\alpha$, one requires that, for some $0< \xi,\eta <1$ (taken to be $\xi=\eta=7/8$ in \cite{FK:mul}) the following series is finite: 
%$$\mathcal{E}\doteqdot \left\{\alpha\in [0,1) :  \suS_{i\notin K_\alpha}\frac{1}{\log^{7/8}q_i}<+\infty\right\},$$
%where $K_\alpha=\{i\in \mathbb{N} :  q_{i+1}< q_i\log^{7/8}q_i\}$, and show that $\lambda(\mathcal{E})=1$ (see Proposition 1.7 in \cite{FK:mul}). This corresponds to Ratner DC with $\xi=\eta=7/8$.
\begin{equation}\label{alphaRatner}
\sum_{k\notin K(\alpha)}\frac{1}{(\log q_n)^\eta} <+\infty, \qquad \text{where}\ K(\alpha):=\{ k\in \mathbb{N}, \ \ a_{k+1}\leq C (\log q_k)^\xi \}.
\end{equation}
In joint work with A.~Kanigowski and J.~Ku{\l}aga-Przymus \cite{KKU:mul}, we were able to generalize this result to higher genus. To do so, it is once again crucial to introduce a suitable Diophantine-like condition, which we called in \cite{KKU:mul} the \emph{Ratner Diophantine Condition} (or RDC) and we describe in \S~\ref{sec:DChigherg}. The main result we prove is the following. 
\begin{theorem}[Kanigowski, Ku{\l}aga-Przymus, Ulcigrai \cite{KKU:mul}] \label{thm:Ratner}
If the restriction of $\varphi_\R \in \mathcal{U}_{\neg min}$ to a minimal component $S_i$ satisfies the Ratner Diophantine Condition, $\varphi_\R: S_i\to S_i$ satisfies the Switchable Ratner Property and is mixing of all orders.
\end{theorem}
\noindent We then show that the RDC is satisfied by almost every IET and therefore can conclude that, for a full measure set of locally Hamiltonian flows in $\mathcal{U}_{\neg min}$, each restriction to a minimal component  is mixing of all orders.

Quantitative estimates on slow, Ratner-type shearing were recently used (in the joint work \cite{KLU:dis} with A.~Kanigowski and M.~Lema{\'{n}}czyk) to study \emph{disjointness  of rescalings}, a property that has recently received a revival of attention in view of its role as possible tool to prove Sarnak Moebius orthogonality conjecture (see the ICM proceedings survey \cite{Le:ICM} and the references therein). In \cite{KLU:dis} we introduce a disjointness criterium based on Ratner shearing and use it (as one of the applications) to show that,  in genus one, typical Arnold flows  have \emph{disjoint rescalings} and satisfy  Moebius orthogonality. Disjointness of rescalings seem to be an important feature of parabolic dynamics: while specific parabolic flows may fail to be disjoint from their rescalings (in primis the horocyle flow on a hyperbolic surface), several recent results seem to indicate  that this property is indeed widespread among parabolic flows (see e.g.~the results in \cite{KLU:dis} on time-changes of horocycle flows). 
In the context of surface flows, disjointness of rescalings has been verified in \cite{BK:dis} for \emph{von Neumann flows} (which can be realized as translation flows on surfaces with boundary). 
Whether one can extend the disjointess result proved in \cite{KLU:dis} for Arnold flows to higher genus smooth flows, remains an open problem and is likely to require a delicate control of Diophantine-like properties. 
%if for any $\kappa\neq 1$ (possibly with finitely many exceptions) $\varphi_\R$ and $\varphi^\kappa_\R$ are \emph{disjoint} in the sense of Furstenberg (see \cite{Le:Moe}). 

%Another feature which seems to emerge as an important feature of parabolic flows 

%Another feature which seems to emerge as an important feature of parabolic flows is  \emph{disjointess  of rescalings}: a rescaling of the flow $\varphi_\R=(\varphi_{t})_{t\in \R}$ is obtained simply by considering  linear time change $\varphi^\kappa_\R:=(\varphi_{\kappa t})_{t\in \R}$.  We say that $\varphi_\R$ has \emph{disjoint rescalings}
%if for any $\kappa\neq 1$ (possibly with finitely many exceptions) $\varphi_\R$ and $\varphi^\kappa_\R$ are \emph{disjoint} in the sen se of Furstenberg (see \cite{Le:Moe}). This property has recently received a revival of attention in view of his role as a possible way to prove Sarnak Moebius orthogonality conjecture (see the survey \cite{Le:Moe} and the references therein). While specific parabolic flows may fail to be disjoint from its rescalings (in primis the horoycle flow on a hyperbolic surface), recent results seem to indicate that this property is indeed widespread among parabolic flows (see e.g.~\cite{FU:tim}). In the context of

\subsubsection*{Polynomial deviations of ergodic averages.}
Slow chaotic behavior manifests itself not only through \emph{slow} mixing, but also through \emph{slow} convergence of ergodic integrals: given an ergodic    area-preserving flow $\varphi_\R$ (or its restriction to an ergodic minimal component $S'\subset S$)  and a real valued observable $f$ with zero-mean, the ergodic integrals $I_T(f,x):=\int_{0}^T f(\varphi_t(x))\ \mathrm{d} t$ decay to zero \emph{polynomially} with some exponent $0<\nu<1$  for almost every initial point, i.e.~
$\big| I_T(f,x) \big| \sim O(T^\nu)$ 
%$\big| \int_0^T f(\varphi_t(p)) \mathrm{ d}\,   t\big| \sim O(T^\nu)$ 
in the sense  that
$\limsup_{T\to \infty} \log {\big|I_T(f,x)\big|}/{\log T}=\nu$. 
%$\limsup_{T\to \infty} {\big|\log \int_0^T f(\varphi_t(p)) \mathrm{ d}\,   t\big|}/{\log T}=\nu$. 
%\begin{equation}\label{nu}
%I_T(f,p)\sim O(T^\nu)\qquad \Leftrightarrow \qquad 
%\limsup_{T\to \infty} {\Big|\log \int_0^T f(\varphi_t(p)) \mathrm{ d}\,   t\Big|}/{\log T}=\nu.
%\lim_{T\to + \infty} I_T(f,p) = \int_M f  \mathrm{ d}\,  \mu, \ \text{whre}\ I_T(f,p)=I_T(f,p,\varphi_\R):=\frac{1} { T} \int_0^T f(\varphi_t(p)) \mathrm{ d}\,   t .
%\end{equation}
This phenomenon, known as \emph{polynomial deviations of ergodic averages}, was discovered experimentally in the $1990s$ by A.~Zorich and explained (for  linear flows on translation surfaces and observables corresponding to cohomology classes) in seminal work by Kontsevitch and Zorich \cite{Zo:dev, Ko:Lya} relating power deviations to Lyapunov exponents of renormalization (see \S~\ref{sec:renorm}). 
Forni in \cite{Fo:dev} could  extend this result to integrals of sufficiently regular functions over translation flows and show that ergodic integrals can display a \emph{power spectrum} of behaviors, i.e.~there are exactly $g$ positive exponents $0< \nu_{g}\leq \cdots \leq \nu_{2}< \nu_1:=1 $ (which correspond to the positive Lyapunov exponents of renormalization)  and, for each, a subspace of finite codimension of smooth observables that present polynomial deviations  as above %in \eqref{nu} 
with exponent $\nu=\nu_i$.  
A finer analysis of the behavior of Birkhoff sums or integrals, beyond the \emph{size} of oscillations, appears in the works \cite{Bu:lim, MMY:aff}: %In \cite{MMY3}, motivated by the study of wandering intervals in affine i.e.m., S.~M., P.~Moussa and J.C.Y. introduced an object called \emph{limit shape} and used it to describe the \emph{shape} of ergodic sums (see \S~3.4 and \S~3.7.3 in  \cite{MMY3}). Roughly speaking these are obtained by looking at suitably rescaled Birkhoff sums, where time is renormalized according to the leading Laypunov exponent of the Kontsevich-Zorich cocycle, whereas the range of the sum  is renormalized using one of the other positive exponents, according to the choice of $f$. After this double rescaling one obtains a sequence of shapes exponentially converging (in the Hausdorff metric)  to the graph of a H\"older function.
 Bufetov   in \cite{Bu:lim} shows in particular that (for  \emph{typical} translation flows and  sufficiently regular observables) the \emph{asymptotic behavior} of ergodic integrals  can be described in terms of $g$ (where $g$ is the genus of the surface) \emph{cocycles} $\Phi_i(t,x)$, $1\leq i\leq g$ (also called \emph{Bufetov functionals}): each $\Phi_i: \R\times S' \to \R$ is a cocycle over the flow $\varphi_\R$ (in the sense that $\Phi_i(t+s,x)= \Phi_i(t,x) + \Phi_i(s,\varphi_t(x))$ for any $ x\in S'$ and $t\in \R$), $\Phi_1(T,x)\equiv T$ and each  $\Phi_i$ has power deviations   $|\Phi_i(T,x)|\sim O(T^{\nu_i})$ with exponent $\nu_i$. Together, the cocycles encode the \emph{asymptotic behavior} of the ergodic integrals up to sub-polynomial behavior, in the sense that, for some constants $c_i=c_i(f)$, 
\begin{equation}\label{eq:expansion}
\int_{0}^T f(\varphi_t(x)) \mathrm{d} t = c_1 T + c_2 \Phi_2(T,x) + \dots + c_g \Phi_g(T,x) + Err(f,T,x),
\end{equation}
 where for almost every $x\in S'$ the \emph{error term} $Err(f,T,p)$ is sub-polynomial, i.e.\ for any $\epsilon>0$ there exists $C_\epsilon>0$ such that $|Err(f,T,p)|\leq C_\epsilon T^\epsilon$. In joint work with Fr{\c a}czek, we recently gave a new proof  of this result  in \cite{FU:Bir} which extends the result to the setting of smooth observables over locally Hamiltonian flows with Morse singularities (in $\mathcal{U}_{min}$ as well as in $\mathcal{U}_{\neg min}$) and also shows   that  the set of locally Hamiltonian flows for which the result holds can be described in terms of a Diophantine-like condition. More precisely, we define in \cite{FU:Bir} the \emph{Uniform Diophantine Condition} (or UDC for short, see \S~\ref{sec:DChigherg}) and show that it has full measure. We then prove the following.
\begin{theorem}[Fr{\c a}czek-Ulcigrai \cite{FU:Bir}]\label{thm:deviations}
If the restriction of the locally Hamiltonian flow $\varphi_\R\in \mathcal{M}$ on a minimal component $S'$ satisfies the \emph{Uniform Diophantine Condition}, for each  $\mathcal{C}^3$ observable $f:S'\to \R$,  there exists $g$ exponents $\nu_i$  and corresponding cocycles $\Phi_i$ such that the expansion \eqref{eq:expansion} holds.
\end{theorem}
%\noindent We return to this condition 

%\subsubsection*{Disjointness of rescalings.}
%Another feature which seems to emerge as an important feature of parabolic flows is  \emph{disjointess  of rescalings}: a rescaling of the flow $\phi_\R=(\varphi_{t})_{t\in \R}$ is obtained simply by considering  linear time change $\varphi^\kappa_\R:=(\varpsi_{\kappa t})_{t\in \R}$.  We say that $\varpsi_\R$ has \emph{disjoint rescalings}
%if for any $\kappa\neq 1$ (possibly with finitely many exceptions) $\varpsi_\R$ and $\varpsi^\kappa_\R$ are \emph{disjoint} in the sense of Furstenberg (see \cite{Le:Moe}). This property has recently received a revival of attention in view of his role as a possible way to prove Sarnak Moebius orthogonality conjecture (see the survey \cite{Le:Moe} and the references therein). While specific parabolic flows may fail to be disjoint from its rescalings (in primis the horoycle flow on a hyperbolic surface), recent results seem to indicate that this property is indeed widespread among parabolic flows (see e.g.~\cite{FU:tim}). In the context of 

% Arnold flows in genus one were also recently shown (by A.~Kanigowski and M.~Lema{\'{n}}czyk and the second author, see \cite{KLU:dis}) to typically have  a property which in particular implies 

\vspace{-1.9mm}
\subsubsection*{Spectral theory.}\label{sec:spectrumlocHam} 
The study of the spectrum of the unitary  operators acting on $L^2(S,\mu)$ given by 
$f\mapsto f\circ \varphi_t$ can shed further light on the chaotic features of the dynamics of the flow $\varphi_\R:= (\varphi_t)_{t\in \R}$ and is at the heart of spectral theory (see \cite{Le:ICM} or \cite{Ul:slo} and the references therein).  
%Spectral theory is the study of the spectrum of the Koopman operator $U_T: L^2()$
While the classification of mixing properties of locally Hamiltonian flows is essentially complete, very little is known on their spectral properties  beyond  the case of genus one (and some sporadic examples, such as \emph{Blokhin examples}, essentially built gluing genus one flows, see  the work  \cite{FL1}). 
The recent result \cite{FFK:Leb}  by   Fayad, Forni and Kanigowski for genus one suggests that it may be possible to prove  that the spectrum is countable Lebesgue also in higher genus when in presence of degenerate, sufficiently strong (multi-saddle) singularities. In the non-degenerate case, though, we recently proved in joint work with Chaika, Fr{\c a}czek and Kanigowski \cite{CFKU:sin} that a \emph{typical} locally Hamiltonian flow on a \emph{genus two} surface with two isomorphic \emph{simple saddles} has {\emph{purely singular}}  {spectrum}. 
% what we believe is the first typical spectral result  for surfaces of higher genus was recently proved in  
This result does not use explicit Diophantine-like conditions, but rather geometry and in particular a special symmetry (the hyperelliptic involution) that surfaces in genus two are endowed with; Liouville-type Diophantine conditions are here imposed by requesting the presence on the surface of large flat cylinders close to the direction of the flow, whose existence for typical flows is then proved by a Borel-Cantelli type of argument (see \cite{CFKU:sin} for details). Extending this result beyond genus two, though, will probably require the use of Rauzy-Veech induction (see \S~\ref{sec:renorm}) and the introduction of new Diophantine-like conditions, which impose some controlled form of degeneration.  The nature of the  spectrum of minimal components of locally Hamiltonian flows in $\mathcal{U}_{\neg min}$ (even in genus one, i.e.~for Arnold flows) is a completely  open problem.
% theory of locally Hamiltonian flows is still largely not understood.
%The first result, to the best of our knowledge, for higher genus surfaces, i.e. for $g\geq 2$, is the following result for genus two (see Figure~\ref{onlysaddles}), which goes in a completely opposite direction.

\subsection{Linearization  and rigidity in higher genus}\label{sec:linhigherg}
A different line of problems in which Diophantine-like conditions in higher genus play a crucial role are conjectures concerning \emph{linearization} and \emph{rigidity} properties of higher genus flows and their Poincar{\'e} sections, GIETs (defined in \S~\ref{sec:higherg}). 
	In analogy with the case of circle diffeos, we say that a GIET $T$ is \emph{linearizable} if it is topologically conjugate to a linear model, namely to a (standard) IET $T_0$.   % \emph{Linear models} for  generalized interval exchange transformations are provided
%, respectively, by (rigid) circle rotations (i.e.~by the map $T_0(x)=x+\alpha \mod 1$ on $I=[0,1]$) and by (standard or classical) \emph{interval exchange transformations} (or for short, IETs).

%\subsubsection*{Linearization and cohomological equation in higher genus}

\subsubsection*{Topological conjugacy and wandering intervals.} To generalize Poincar{\'e} and Denjoy work, one needs first of all a combinatorial invariant which extends the notion of rotation number. Such an invariant can be constructed by recording the combinatorial data of a renormalization process, as we explain in \S~\ref{sec:renorm}. One of the crucial differences between GIETs and circle diffeomorphisms, though, is the failure  of a generalization of Denjoy theorem: there are smooth GIETs that are semi-conjugate to a minimal IET for which the semi-conjugacy is \textit{not} a conjugacy, in other words they have wandering intervals (see the examples found in \cite{CG:aff, BHM:per} in the class of periodic-type (affine) IETs and, more generally, \cite{MMY:aff}). 
%a  phenomenon, first discovered by Levitt but a a non-uniquely ergodic example \cite{Le}, then explored  by Camelier Gutierrez \cite{ CamelierGutierrez}, Cobo \cite{Cobo} and Bressaud, Hubert and Mass \cite{BressaudHubertMaass} for %special (families of) uniquely ergodic 
%examples of periodic type in the sense of \S~\ref{sec:DChigherg} and finally by \cite{MMY2}). 
It is important to stress that this is  \emph{not} a low-regularity phenomenon, nor it is related to special  arithmetic assumptions: as shown by the key  work  \cite{MMY:aff} by Marmi, Moussa and Yoccoz,  wandering intervals exist
 even for piecewise \emph{affine} (hence analytic) GIETs (called AIETs), for almost every topological conjugacy class. The presence of wandering intervals is on the contrary  expected to be \emph{typical} (see e.g.~the conjectures in \cite{MMY:lin, Gh:une})  and it is closely interknit with the absence of a \emph{Denjoy Koksma inequality} and more in general \emph{a priori bounds} for renormalization, see~\cite{GU:rig}.
% for AIETs with \emph{almost every} rotation number, 
%as shown later in the  work \cite{MMY2} by Marmi, Moussa and Yoccoz, which actually indicates that the existence of wandering interval is in some sense \emph{typical}\footnote{See for example the statement of Proposition~\ref{AIETwi}, which is taken from \cite{MMY}: wandering intervals are shown to exist for a full measure set of rotation numbers as long as the log-slope vector of the AIET has a typical projection on the Oseledets filtration (and, conjecturally, as long as it projects on any positive Oseledets exponent).}.
%the absence of a \emph{Denjoy Koksma inequality}\footnote{We recall that the \emph{Denjoy-Koksma inequality} is an ergodic-theoretic statement which gives boundedness of Birkhoff sums of bounded variation observables at special times: given $f:I\to I$ is a function of bounded variation on $I=[0,1]$ and $R_\rho$ a rotation by $\rho$, if $p_n/q_n$ are the \emph{convergents} of $\rho$ (given by $p_n/q_n=[a_0,\dots, a_n]$ where $[a_0,\dots, a_n,\dots]$ is the continued fraction expansion of $\rho$), the Birkhoff sums $\sum_{k=0}^{q_n}f(x)$ at times $q_n$ are uniformly bounded, independently on $n\in\mathbb{N}$ and $x\in I$.} and \emph{a priori bounds}. This has far-reaching consequences, the most spectacular of which being 

\subsubsection*{Local obstructions to  linearization.} 
As an important first step towards local linearization, we already mentioned the \emph{cohomological equation} $\varphi\circ T- \varphi=\phi$ in \S~\ref{sec:g1}, where $T=R_\alpha$ was a rotation.  Whether the cohomological equation could be solved when $T$ is a IET, under suitable assumptions, was unknown until the pioneering work of Forni \cite{Fo:coh}, % (see also \cite{Forni2}), 
who brought to light the existence of a  \textit{finite} number of obstructions to the existence of a (piecewise finite differentiable) solution. We remark that obstructions to solve the cohomological equation 
%(previously known for skew products over rotations, see \cite{Ka:com}),  
have been since then discovered to be a characteristic phenomenon in \emph{parabolic dynamics} (e.g.~their existence have been proved by Flaminio and Forni for horocycle flows \cite{FF:hor} and nilflows on nilmanifolds \cite{FF:nil}, 
%, see for example the works by Flaminio and Forni on the cohomological equation for horocycle flows    and for  nilflows on nilmanifolds or. see 
see also the ICM talk~\cite{Fo:asy}). 
Forni's work is a breakthrough that paved the way for the development of a linearization theory in higher genus. 
 
Another  breakthrough, which put the stress on the \emph{arithmetic} aspect of linearization in higher genus,  was  achieved by  Marmi-Moussa-Yoccoz in their work \cite{MMY:lin} (and related works \cite{MMY:coh, MY:Hol}). 
 In \cite{MMY:coh}, in particular, they reproved and extended Forni's result using the IETs renormalization described in \S~\ref{sec:DChigherg} and introduced the \emph{Roth-type} condition (see also \S~\ref{sec:DChigherg}),  as an explicit Diophantine-like condition on the IET needed to solve the cohomological equation $\varphi\circ T- \varphi=\phi-\xi$, where $\xi$ is a piecewise constant function which embodies the finite dimensional \emph{obstructions}.   
 % (after removing a \emph{correction} living in a finite dimensional subspace).  
% they reproved Forni's result on the cohomological equation using a renormalization approach based on Rauzy-Veech induction, thus describing explicitely a full measure Diophantine-like condition on the IET (a condition that, in analogy with rotations, they called \emph{Roth type}, see \cite{MMY} and also \cite{MarmiYoccoz} for a variation of this condition). 
%A refinement of this result by Marmi, Moussa and Yoccoz (with an explicit full measure arithmetic condition on $T_0$ called \emph{Roth-type}) 
This result, combined with a generalization of  Herman's \emph{Schwarzian derivative trick},  then led to the proof in \cite{MMY:lin} by the same authors that, for any $r\geq 2$, the $\mathcal{C}^r$ \emph{local conjugacy class} of almost every IET $T$ (more precisely, of  any $T$ of \emph{restricted} Roth-type, see \S~\ref{sec:DChigherg})  is a \emph{submanifold of finite codimension}. Marmi, Moussa and Yoccoz also conjectured that for $r=1$ it is a submanifold of codimension $(d-1)+(g-1)$, where $d$ is the number of exchanged intervals and $g$ the genus of the surface of which $T$ is a Poincar{\'e} section. For the measure zero class of IETs of \emph{hyperbolic periodic type} (see \S~\ref{sec:DChigherg}), this conjecture has recently been proved by  Ghazouani in \cite{Gh:loc}. The proof of this result for almost every IET will require the introduction of a new suitable Diophantine-like condition on IETs.

\subsubsection*{Rigidity of GIETs.} We say that a class of (dynamical) systems is \emph{geometrically rigid} (or also $\mathcal{C}^1$-rigid), if the existence of a topological conjugacy between two objects in the class automatically imply that the conjugacy is actually $\mathcal{C}^1$. The global linearization results by Herman and Yoccoz recalled in \S~\ref{sec:g1} shows that the class of (smooth, or at least $\mathcal{C}^3$) circle diffeomorphisms with Diophantine rotation number is geometrically rigid (and actually $\mathcal{C}^\infty$-rigid, i.e.~if a smooth circle diffeo is conjugated via a homeomorphism $h$ to $R_\alpha$ with $\alpha $ satisfying the DC, then $h$ is $\mathcal{C}^\infty$).  %: if a circle diffeo is topologically conjugated to a linear model $R_\alpha$, where $\alpha$ is its rotation nu
 We already saw that this can be reinterpreted as a rigidity result for flows on surfaces of genus one (see Theorem~\ref{HYthm}). In joint work with S.~Ghazouani, we recently proved a generalization of this result to genus two.
\begin{theorem}[Ghazouani, Ulcigrai \cite{GU:rig}]
Under a full measure  Diophantine-like condition, a foliation on a  \emph{genus two} surface which  is topologically conjugate to the foliation given by a linear flow with Morse saddles, is also $\mathcal{C}^1$ conjugate to it.
\end{theorem}
\noindent Here full measure refers to the Katok measure class on the linear flow models (see the definition given earlier in this section). For simplicity, we stated the result for flows with simple, Morse-type saddles; degenerate saddles can also be considered, but then one has to further assume that the foliations are \emph{locally} $\mathcal{C}^1$ conjugated in a neighborhood of the multi-saddle. 
Both these results can be reformulated at level of Poincar{\'e} sections: we introduce more precisely a rather subtle Diophantine-like conditions on (irreducible) IETs of any number of intervals $d\geq 2$, that we call the\emph{ Uniform Diophantine Condition}, or UDC (we comment on it in \S~\ref{sec:higherg}) and show that it is satisfied by almost every (irreducible) IET on $d$. We then prove:
\begin{theorem}[Ghazouani, Ulcigrai \cite{GU:rig}]\label{rigidityg2}
 If an irreducible $d$-IET  $T_0$ with $d=4$ or $d=5$ satisfies the \text{UDC},  then any  $\mathcal{C}^3$-generalized interval exchange map $T$ which is topologically conjugate to $T_0 $, and whose  \emph{boundary} ${B}(T)$ \emph{vanishes},   is actually conjugated to $T_0$ via a $\mathcal{C}^1$ diffeomorphism.
% then the conjugacy between $T$ and $T_0$ is actually a  of $[0,1]$ of class  
 %In other words, almost every standard irreducible IET with $4$ or $5$ continuity intervals is  geometrically rigid.
% $T_0$  is also $\mathcal{C}^1$ conjugate to it.
\end{theorem}
\noindent The \emph{boundary operator} ${B}(T)$ which appears in this statement is a $\mathcal{C}^1$-conjugacy invariant introduced in \cite{MMY:lin}; it encodes the holonomy at singular points of the surface of which $T$ is a Poincar{\'e} section. Requesting that ${B}(T)$ vanishes is therefore a necessary condition for the existence of a conjugacy of class $\mathcal{C}^1$. 
Theorem~\ref{rigidityg2} solves for $d=4,5$ one of the open problems suggested  by Marmi, Moussa and Yoccoz in \cite{MMY:lin}, where they conjecture the result to hold also for any other larger $d$. The result which is missing to prove the conjecture in its generality is a generalization of an estimate used in \cite{MMY:aff} to show existence of wandering intervals in affine IETs. The main result in \cite{GU:rig}, on the other hand (namely a dynamical dichotomy for the orbit of $T$ under renormalization) is already proved for IETs which satisfy the UDC for any $d\geq 2$.

\section{Renormalization and cocycles}\label{sec:renorm}
In this section we introduce the renormalization dynamics which is used as main tool to impose   Diophantine-like  conditions in higher genus.
%we will exploit in the next \S~\ref{sec:DChigherg} to describe the Diophantine-like  conditions behind some of the results which we summarized in the previous \S~\ref{sec:higherg}. 
%are  expressed exploiting a \emph{renormalization} algorithm, in order to describe them (in ) we first introduce in this section the idea of renormalization in this setting and reinterpret classical Diophantine conditions for rotation numbers from this perspective. type of conditions are imposed, give some examples and comment on their nature. Since these conditions are expressed exploiting a \emph{renormalization} algorithm, we first introduce the idea of renormalization in this setting and reinterpret classical Diophantine conditions for rotation numbers from this perspective.
Renormalization in dynamics is a powerful tool to study dynamical systems which present forms of self-similarity (exact or approximate) at different scales. A map $T:I\to I$ of the unit interval which is (infinitely)  \emph{renormalizable} is such that one can find a (infinite) sequence of nested subintervals $I_{n+1}\subset I_n\subset \dots \subset I$ such that the \emph{induced dynamics} $T_n: I_n\to I_n$ (obtained by considering the first return map of $T$ on $I_n$) is well defined and, up to \emph{rescaling}, belongs to the same class of dynamical systems of the original $T$. Here, the rescaling, which is done so that the \emph{rescaled} (or \emph{renormalized}) map  acts again on an interval of unit length,  is given by the map $x\mapsto T_n(|I_n|x)/|I_n| $. We will now describe renormalization in the context of rotations first and then IETs. In both cases, at the level of (minimal) flows (or equivalently orientable foliations) on surfaces, the inducing process corresponds to taking shorter and shorter Poincar{\'e} sections of a given surface flow (on the torus or on a higher genus surface).
%We will see now that irrational rotations (and more in general circle diffeomorphisms with irrational rotation number)  and interval exchange maps (and more in general GIETs with \emph{irrational} rotation number, see \S~\ref{sec:rigidity}) are (infinitely) renormalizable.

\subsubsection*{Renormalization algorithms.} 
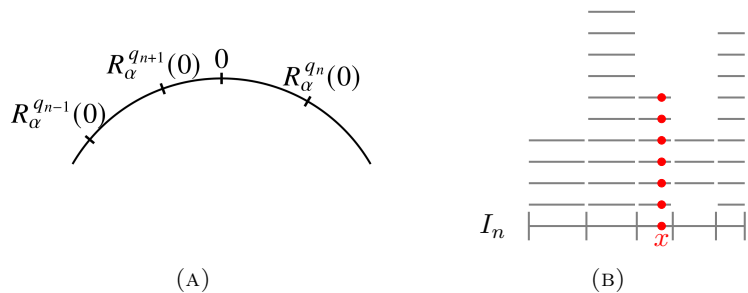
\begin{figure}[h!]
\subfloat[\label{RenRotation}]{\includegraphics[width=5cm]{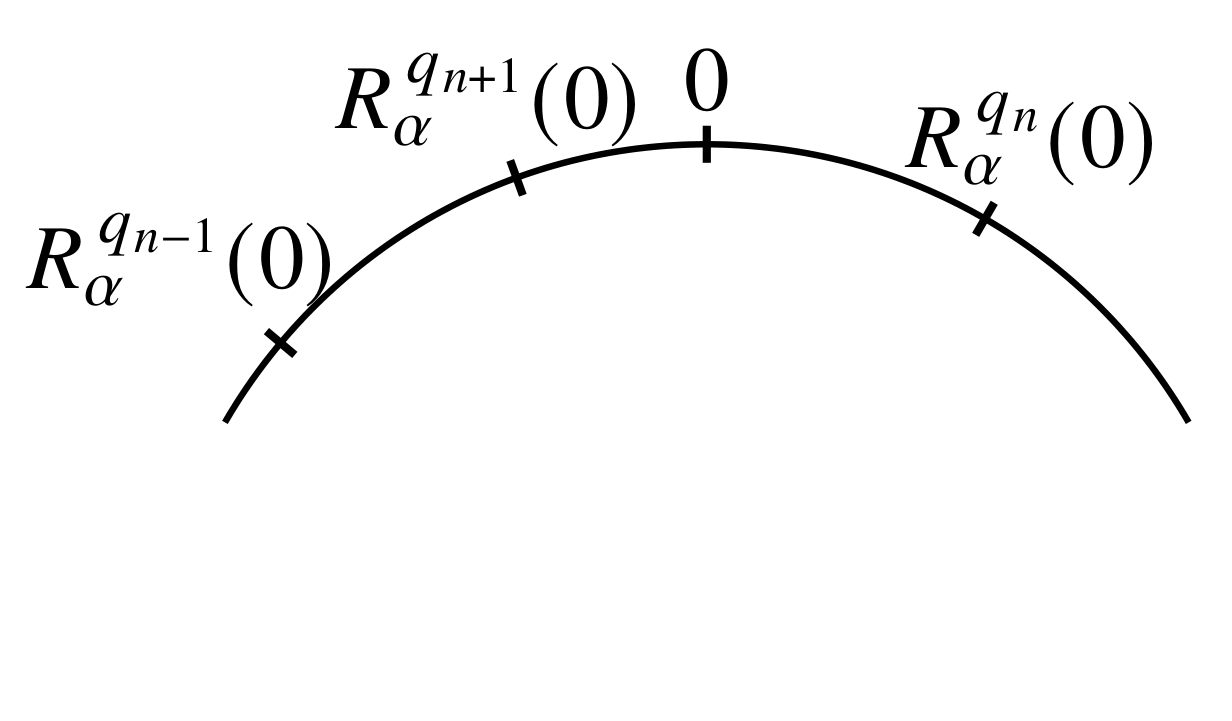}}
%\centering
%\resizebox{160pt}{!}{ 
%image size
%\begin{tikzpicture}
%\filldraw [fill=white] (0,0) circle [radius=6cm]; %draws full circle
%\draw[black, line width= 2pt] (6,0) arc (0:180:6); %thickness of circle rc
%\draw[white, line width= 5pt] (6,0) arc (0:-210:6); %length of circle rc n the left
%\draw[white, line width= 5pt] (6,0) arc (0:30:6); %length of circle arc n the right
%\draw[line width=2.5pt] (90:5.8cm) -- (90:6.2cm);
%\node[font=\Huge] at (90:6.7cm) {$0$};
%\draw[line width=2.5pt] (60:5.8cm) -- (60:6.2cm);
%\node[font=\Huge] at (60:x7cm) {$R_\alpha^{\, q_n}(0)$};
%\draw[line width=2.5pt] (110:5.8cm) -- (110:6.2cm);
%\node[font=\Huge] at (110:6.9cm) {$R_\alpha^{\, q_{n+1}}(0)$};
%\draw[line width=2.5pt] (140:5.8cm) -- (140:6.2cm);
%\node[font=\Huge] at (140:7.4cm) {$R_\alpha^{\, q_{n-1}}(0)$};
%\end{tikzpicture}
%}}
\hspace{9mm}
\subfloat[\label{RohlinTower}]{%\includegraphics[width=6.6cm
\resizebox{110pt}{!}{ %image size
\begin{tikzpicture}
\node at (3,0) {$I_n$}; %Bottom line
\draw[gray, thick] (3.5,0) -- (6.5,0);
\draw[gray, thick] (3.5,-0.2) -- (3.5, 0.2);
\draw[gray, thick] (4.3,-0.2) -- (4.3, 0.2);
\draw[gray,  thick] (5,-0.2) -- (5, 0.2);
\draw[gray, thick] (5.5,-0.2) -- (5.5, 0.2);
\draw[gray, thick] (6.1,-0.2) -- (6.1, 0.2);
\draw[gray, thick] (6.5,-0.2) -- (6.5, 0.2);
\draw[gray, thick] (3.5, 0.3) -- (4.275, 0.3); %first tower
\draw[gray, thick] (3.5, 0.6) -- (4.275, 0.6); 
\draw[gray, thick] (3.5, 0.9) -- (4.275, 0.9); 
\draw[gray, thick] (3.5, 1.2) -- (4.275, 1.2); 
\draw[gray, thick] (4.325, 0.3) -- (4.975, 0.3); %second tower
\draw[gray, thick] (4.325, 0.6) -- (4.975, 0.6); 
\draw[gray, thick] (4.325, 0.9) -- (4.975, 0.9); 
\draw[gray, thick] (4.325, 1.2) -- (4.975, 1.2);
\draw[gray, thick] (4.325, 1.5) -- (4.975, 1.5);
\draw[gray, thick] (4.325, 1.8) -- (4.975, 1.8); 
\draw[gray, thick] (4.325, 2.1) -- (4.975, 2.1); 
\draw[gray, thick] (4.325, 2.4) -- (4.975, 2.4);
\draw[gray, thick] (4.325, 2.7) -- (4.975, 2.7); 
\draw[gray, thick] (4.325, 3) -- (4.975, 3);
\draw[gray, thick] (5.025, 0) -- (5.475, 0); %third tower
\draw[gray, thick] (5.025, 0.3) -- (5.475, 0.3);
\draw[gray, thick] (5.025, 0.6) -- (5.475, 0.6); 
\draw[gray, thick] (5.025, 0.9) -- (5.475, 0.9); 
\draw[gray, thick] (5.025, 1.2) -- (5.475, 1.2);
\draw[gray, thick] (5.025, 1.5) -- (5.475, 1.5);
\draw[gray, thick] (5.025, 1.8) -- (5.475, 1.8);%\draw[gray, very thick] (5.525, 0.3) -- (6.075, 0.3); %fourth tower
\draw[gray, thick] (5.525, 0.6) -- (6.075, 0.6); 
\draw[gray, thick] (5.525, 0.9) -- (6.075, 0.9); 
\draw[gray, thick] (5.525, 1.2) -- (6.075, 1.2);

\draw[gray, thick] (6.125, 0.3) -- (6.5, 0.3); % sixth tower
\draw[gray, thick] (6.125, 0.6) -- (6.5, 0.6); 
\draw[gray, thick] (6.125, 0.9) -- (6.5, 0.9); 
\draw[gray, thick] (6.125, 1.2) -- (6.5, 1.2);
\draw[gray, thick] (6.125, 1.5) -- (6.5, 1.5);
\draw[gray, thick] (6.125, 1.8) -- (6.5, 1.8); 
\draw[gray, thick] (6.125, 2.1) -- (6.5, 2.1); 
\draw[gray, thick] (6.125, 2.4) -- (6.5, 2.4);
\draw[gray, thick] (6.125, 2.7) -- (6.5, 2.7);
\filldraw[red] (5.345,0) circle (1.5pt); %x and points on levels ofiddle tower
\node[red] at (5.345,-0.2) {$x$};
\filldraw[red] (5.345,0.3) circle (1.5pt);
\filldraw[red] (5.345,0.6) circle (1.5pt);
\filldraw[red] (5.345,0.9) circle (1.5pt);
\filldraw[red] (5.345,1.2) circle (1.5pt);
\filldraw[red] (5.345,1.5) circle (1.5pt);
\filldraw[red] (5.345,1.8) circle (1.5pt);
\end{tikzpicture}
}}
 \caption{Renormalization algorithms for rotations and IETs.\label{renormalization}}
\end{figure}

If $T=R_\alpha$ is a rotation by an irrational $\alpha$ and $q_n$, $n\in \mathbb{N}$, are the denominators of the convergents $p_n/q_n$ of $\alpha$, then one can consider as sequence $(I_n)_{ n\in \mathbb{N}}$ the  shrinking arcs on $S^1$ which have as  endpoints  $R_{\alpha}^{q_n}(0)$ and $R_{\alpha}^{q_{n+1}}(0)$, that correspond dynamically to consecutive closest returns of the orbit of $0$ (see Figure~\ref{RenRotation}). The induced map $T_n$ is then again a rotation $R_{\alpha_n}$, with rotation number $\alpha_n=\mathcal{G}^n(\alpha)$, where  $\mathcal{G}$ is the Gauss map $\mathcal{G}(x)=\{1/x\}$. 
%Let us recall that $\mathcal{G}$ admits a finite invariant measure (the Gauss measure). 

 Similarly, for a $d$-IET $T$, one wants to choose the nested sequence $(I_n)_{n\in \mathbb{N}}$ of inducing intervals so that the induced maps $T_n$ are all IETs of the same number $d$ of subintervals. Given any minimal $T$ (or more generally any IET satisfying the \emph{Keane condition}, i.e.~such that the orbits of its discontinuity points are infinite and distinct), classical 
%Renormalization operators for IETs $T:I\to I$ on $I=[0,1]$, %in $\mathcal{X}^r $ 
%produce a choice of subintervals $I^n$ such that $T_n$ is an IET of the same number $d$ of subintervals. 
%of the same number of intervals.  The image $\mathcal{R}(T)$ of $T$ under the renormalization operator is
% then by definition the GIET acting on $I=[0,1]$ obtained by \emph{normalising}, i.~.e.~conjugating by the affine transformation which maps $I'$ to $I$, so that the image is again a GIET on $I$. 
algorithms which produce such an infinite sequence ${(I_n)}_{n\in \mathbb{N}}$ 
are %\emph{renormalize} IETs is %the \emph{Rauzy-Veech algorithm}, also called 
the \emph{Rauzy-Veech induction algorithm} % (whose definition we recall in \S~\ref{Rauzysec}), 
 (see Veech \cite{Ve:gau} or \cite{Yoc:Clay}  %and used starting from the seminal papers by 
and the references therein) and \emph{Zorich induction}, an acceleration of the same algorithm introduced by Zorich in \cite{Zo:gau}.  For the definitions of these algorithms, which we will not use in the following, we refer the interested reader to the lecture notes \cite{Yoc:Clay}. 
% (see e.g.~\cite{Zo:dev,  AF:wea, Ch:dis}). 
%Inducing $T$ to $I^1$ and then renormalizing $T^1$, one can an operator 
% order for the associated renormalization operator
% $\mathcal{R}: \mathcal{I}^d\to  \mathcal{I}^d$ on the space of $d$-IETs. For $\mathcal{R}$ to preserve 
% 
One can show that, for $d=2$, Zorich induction  corresponds  to the renormalization of rotations given by the  Gauss map.

On the parameter space  $\mathcal{I}_d$ of all $d$-IETs, these algorithms induce   renormalization operators $\mathcal{R}: \mathcal{I}_d \to \mathcal{I}_d$, which associate to $T$ the $d$-IET $\mathcal{R}(T)$ obtained by applying one step of the corresponding induction and then renormalizing the induced map to act on $[0,1]$. Veech showed that Rauzy-Veech renormalization admits a conservative absolutely continuous invariant measure, that induces a  \emph{finite} invariant measure for the \emph{Zorich acceleration}, as proved in  \cite{Zo:gau}.  The ergodic properties of the renormalization dynamics in parameter space  have been intensively studied and are by now well understood, see e.g.~%\cite{Zo:gau,  Bu:dec, AGY, AB:exp}, or
 \cite{Yoc:B} and the references therein for a brief survey.
 
%if the permutation is irreducible and the so called \emph{Keane condition}, which guarantees that orbits of singularities are disjoint, is satisfied), the renormalization can be  

\subsubsection*{Rohlin towers and matrices} 
After $n$ steps of induction, one can recover the original dynamics through the notion of \emph{Rohlin towers} as follows: 
if $I_n^i$ is one of the subintervals of $T_n$ and $r:=r_n^i$ is its first return time to $I_n$ under the action of $T$, the intervals $I_n^i, T(I_n^i), \dots, T^{r-1}(I_n^i)$ are disjoint. Their union is called a \emph{Rohlin tower} of step $n$ and each of them is called a \emph{floor} (see Figure~\ref{RohlinTower}). Given an infinitely renormalizable $T$, for any $n$ one can see $[0,1]$ as a union of $d$ Rohlin towers of step $n$, as shown in Figure~\ref{RohlinTower}. Rohlin towers thus produce  a sequence of \emph{partitions} of $[0,1]$ (into floors of towers of step $n$).

Renormalization produces also a sequence of $d\times d$ matrices $A_n$, $n\in\mathbb{N}$, with integer entries, which should be thought of as  \emph{multi-dimensional continued fraction} digits and describe \emph{intersection numbers} of Rohlin towers.  The matrices $(A_n)_{n\in\N}$ are defined so that the entries  of the product $A^n:= A_n\cdots A_1$ have the following dynamical meaning: the $(i,j)$ entry $(A^n)_{ij}$ is the number of visits of the orbit of any point $x\in I_n^j$ to the initial subinterval $I_0^i$ until its first return time $r_n^j$; % $I^{n}_i$ until , namely of $\{ T^\ell I^n_i, 0\leq \ell <n\} $ o to , or,
in other words, $(A^n)_{ij}$  is the number of floors of the $j^{th}$ tower of level $n$ which are contained in $I_0^i$.  These entries generalize the classical continued fraction digits: for $d=2$, indeed, the matrices  $(A_n)_{n\in \mathbb{N}}$ associated to $R_\alpha$, for $n$ of alternate parity, have respectively the form 
$$
\begin{pmatrix} 1 & a_n \\ 0 & 1 \end{pmatrix} \quad \text{or} \quad \begin{pmatrix} 1& 0 \\ a_n & 1 \end{pmatrix},
$$
where $a_n$ are the entries of the continued fraction expansion $\alpha=[a_0,a_1,\dots, a_n ,\dots]$. 
%Similarly to arithmetic conditions for real numbers, 
Diophantine-like conditions for IETs are  defined by imposing conditions on these matrices, on their growth as  well as on their hyperbolicity,  %we will explain 
see in \S~\ref{sec:DChigherg}. The matrices $(A_n)_{n\in \mathbb{N}}$ are produced by the renormalization dynamics: for rotations,  the entries $(a_n)_{n\in\N}$ of the continued fraction  expansion of $\alpha$ satisfy $a_n=a_0(\mathcal{G}^n(\alpha))$, where $a_0(\cdot )$ is an integer valued function on $[0,1]$. Similarly, one has now that $A_n=A_0(\mathcal{R}^n (T))$, where $A_0: \mathcal{I}_d\to SL(d,\mathbb{Z})$ is a matrix valued function on the space $\mathcal{I}_d$ of $d$-IETs, i.e.\ a \emph{cocycle} (known as the \emph{Rauzy-Veech cocycle},  or \emph{Zorich cocycle} if considering the Zorich acceleration).

\subsubsection*{Positive and balanced accelerations.}
It turns out though, that Zorich acceleration is often not sufficient (see for example \cite{KM:bt} and \cite{Ki:typ} where it is shown that the classical Diophantine notions of bounded-type \cite{KM:bt} and Diophantine-type \cite{Ki:typ} do not generalize naturally when using Zorich acceleration).  
Two accelerations which play a key role in Diophantine-like conditions are the \emph{positive} and the \emph{balanced acceleration}. By \emph{accelerations} we mean here an induction which is obtained by considering only a subsequence $(n_k)_{k\in \N}$ of Rauzy-Veech times. The associated (accelerated) cocycle is then obtained considering products 
$$A(n_{k}, n_{k+1}):= A_{n_{k+1}-1}\cdots A_{n_k+1} A_{n_{k}}.$$
The \emph{positive acceleration} appears in the works by Marmi, Moussa and Yoccoz \cite{MMY:coh, MMY:lin, MY:Hol}. They showed that if $T$ satisfies the Keane condition, for any $n$ there exists $m>n$ such that $A(n,m)$ is a strictly positive matrix. The accelerated algorithm then corresponds to choosing the sequence $(n_k)_{k\in \mathbb{N}}$ setting $n_0:=0$ and then, for $k\geq 1$,  choosing $n_k$ to be the smallest integer $n>n_{k-1}$ such that $A(n_{k-1}, n)$ is strictly positive. On the other hand, to define the \emph{balanced acceleration}, one considers a subsequence $(n_k)_{k\in\N}$ of Rauzy-Veech times $n $ for which the corresponding Rohlin towers are \emph{balanced}, in the sense that ratios of widths $|I_{n}^i|/|I_{n}^j|$ and heights $r^{n}_i/r^n_j$ are uniformly bounded above and below. We will return to these accelerations and some instances in which they are helpful in \S~\ref{sec:DChigherg}.

\subsubsection*{Combinatorial rotation numbers.}
We remark that the definition of Rauzy-Veech induction can be extended also to a GIET $T$ (under the Keane condition, which guarantees that $\mathcal{R}^n(T)$ can be defined for every $n\in\N$) and then  exploited to give a combinatorial notion of \emph{rotation number} as well as a definition of \emph{irrationality} in higher genus  (following \cite{MMY:lin, MY:Hol}, see also \cite{Yoc:Clay}).  As one computes the induced maps $(T_n)_{n\in \N}$, one can indeed record the sequence $(\pi_n)_{n\in \N}$ of permutations of the GIETs $(T_n)_{n\in\N}$: this  sequence provides the desired \emph{combinatorial rotation number} for $d>2$. We say that a GIET is \emph{irrational} if the sequence of matrices $(A_n)_{n\in \N}$ have a positive acceleration (or equivalently, in the terminology introduced by Marmi, Moussa and Yoccoz, the path described by $(\pi_n)_{n\in \N}$ is \emph{infinitely complete}). One can then show that two \emph{irrational} GIETs with the same rotation number are \emph{semi-conjugated} (see e.g.\ \cite{Yoc:Clay}), a result that generalizes a property of rotations numbers and circle diffeos and hence explains the choice of calling this higher genus combinatorial object the `\emph{rotation number}' of a GIET

.
%that we call Poincar{\'e}-Yoccoz theorem (see Theorem~\ref{thm:PY}). 
%Renormalization provides also  combinatorial i 

\subsubsection*{Renormalization of Birkhoff sums}
Given $T:I\to I $ and a function $f:I\to \R$, we denote by $S_n f:= \sum_{k=0}^{n-1} f\circ T^k$ the $n^{th}$-Birkhoff sum (of the function $f$ under the action of $T$). When $T=R_\alpha$ is a rotation (or a circle diffeo), it is standard to study first Birkhoff sums of the form $S_{q_n} f$ for $q_n$  convergent of $\alpha$, corresponding to closest returns, and then use them to \emph{decompose} more general Birkhoff sums. Similarly, renormalization for (G)IETs can be exploited to produce \emph{special Birkhoff sums}, namely Birkhoff sums of a special form that can be understood first, exploiting renormalization, and then used to decompose and study general Birkhoff sums.  For each $n\in \N$, if $T_n:I_n\to I_n$ is the induced map after $n$ steps of renormalization, the $n^{th}$ \emph{special Birkhoff sum} is the induced function $S(n)f: I_n\to I_n$, defined by $S(n)f(x):= S_{r_n^i} f (x)$ if $x\in I_n^i$. Thus, since $r_n^i$ is the height of the Rohlin tower over $I_n^i$, 
 the value  $S(n)f(x)$ is obtained summing the orbit \emph{along the tower}  which has $x$ in the base, see Figure~\ref{RohlinTower}. Notice that for $d=2$, when considering Zorich acceleration, these reduce to sums of the form $S_{q_n} f(x)$. 
The associated \emph{special Birkhoff sums operators} $S(n)$, $n\in\N$, map $f: I\to \R $ to $S(n)f: I_n\to \R$. When $f$ is piecewise constant and takes a constant value $f_i$ on each $I_i$, $S(n)$ can be identified with a linear operator given by the (studied  acceleration of the) Rauzy-Veech cocycle $A^n = A_{n}\cdots A_1$ as follows: one can show that  $S(n)f$ takes constant values $f^n_i$ on each $I_n^i$ and the column vectors $\underline{f}:=(f_i)_{i=1}^d$ and $\underline{f}^n:=(f^n_i)_{i=1}^d$ are related by $\underline{f}^n = A^n \, \underline{f}$. Thus, special Birkhoff sums operators can be seen as infinite dimensional extensions of the Rauzy-Veech cocycle (and its accelerations).

When considering a rotation $R_\alpha$, to decompose $S_n f(x)$ into Birkhoff sums of the form $S_{q_k} f(x_j)$, one can write $n=\sum_{k=0}^{k_n} b_k q_k$, where $k_n$ is the largest integer $k$ such that $q_{k}<n$ and $b_k$ are integers such that $0\leq b_k\leq a_{k_1}$ (a factorization sometimes known as \emph{Ostrowsky decomposition}). Correspondingly, recalling that $S_{q_k} f(x_j)= S(k) f (x_j)$ when $x_{j}^k\in I_k$, we can write 
\begin{equation}\label{decomposition}
 S_n f(x)= \sum_{k=0}^{k_n}\sum_{j=0}^{b_k-1} S(k) f(x^k_j) , \qquad \text{where} \ x_{j}^k\in I_k, \quad \text{for\ all}\ 0\leq j<b_k. 
\end{equation}
For IETs one can also get an analogous decomposition of any Birkhoff sums $S_n f(x)$ into special Birkhoff sums, which has the same form \eqref{decomposition}, but where $0\leq b_k\leq \Vert A^n\Vert:=\sum_{i,j} (A^n)_{ij}$ and the decomposition is obtained \emph{dynamically}, by decomposing the orbit of $x$ until time $n$ into blocks, each of which is contained in a tower and hence corresponds to a special Birkhoff sums.

\subsubsection*{Renormalization in moduli spaces.}
We conclude this section mentioning that these renormalization algorithms (for rotations and IETs) describe a discretization of a renormalization dynamics on the moduli space of surfaces. In genus one, the Gauss map is well known to be related to the geodesic flow on the modular surface (which can be seen as the moduli space of flat tori), see e.g.~\cite{Se:mod}. Similarly, (an extension of) Rauzy-Veech induction can be obtained as Poincar{\'e} map of the \emph{Teichmueller geodesic flow} on the moduli space of translation surfaces (see e.g.~\cite{Zo:fla}). 

The full measure Diophantine-like conditions that we discuss in this survey are satisfied by  (Poincar{\'e} maps of) linear flows in almost every direction on almost every translation surface in these moduli space (with respect to the Lebesgue, or Masur-Veech measure, see \cite{CW:ICM}). A different question is whether these properties hold for a \emph{given} surface in almost every direction, in particular if the surface has special properties, for example is a torus cover (i.e.~it is a square-tiled surface), or has special symmetries (e.g.~it is a \emph{Veech surface} or it belongs to a n $SL(2,\R)$-invariant locus, see \cite{CW:ICM}). In these settings, while some results can be obtained by general measure-rigidity techniques (in particular from the work \cite{EC:Ose} by Chaika and Eskin, see also the ICM proceedings \cite{CW:ICM} and the references therein),  to describe explicit Diophantine-like conditions 
 it is often helpful to exploit or develop \emph{ad-hoc} renormalization algorithms (for example one can use finite extensions of the Gauss map to study square-tiled surfaces, see e.g.~\cite{MMY:sqt}, or construct Gauss-like maps for some Veech surfaces, see e.g.~\cite{SU:oCF}). 

%A full measure set of directions can be sometimes obtained through results such as \cite{EC

%These moduli spaces contain also many natural locally affine submanifolds (also called \emph{loci}, invariant under the natural $SL(2,\R)$-action, see \cite{CW:ICM} and the references therein), the simples of which are generated by Veech surfaces. It is often interesting to 

\section{Diophantine-like conditions in higher genus}\label{sec:DChigherg}
We finally describe in this section some of the Diophantine-like conditions which were introduced to prove some of the results on typical ergodic and spectral properties of smooth area-preserving flows on surfaces (see \S~\ref{sec:locHam}) and on \emph{linearization} (such as solvability of the cohomological equation and  rigidity questions in higher genus, see \S~\ref{sec:linhigherg}). 
%require to impose \emph{typical} conditions on the (G)IETs which arise as Poincar{\'e} map, often in the form of \emph{Diophantine-like} conditions. 

\subsection{Bounded-type IETs and Lagrange spectra}\label{sec:bt}
We start with two important classes of IETs, namely \emph{periodic-type} and \emph{bounded-type} IETs,  both of which have measure  zero  in the space $\mathcal{I}_d$ of IETs (although full Hausdorff dimension in the case of bounded-type IETs), but often constitute an important class of IETs in which dynamical and ergodic properties can be tested.  

%which Bounded-type $d$-IETs have measure zero but nevertheless 
%correspond to periodic-type and bounded-type in genus one.
One of the simplest requests on a (G)IET is that its orbit under renormalization is \emph{periodic}, so that the sequence of Rauzy-Veech cocycle matrices $(A_n)_{n\in \mathbb{N}}$ introduced in the previous section \S~\ref{sec:renorm} is \emph{periodic}, i.e.\ there exists $p>0$ such that $A_{n+p}=A_n$ for every $n\in \mathbb{N}$. We will furthermore request that the \emph{period matrix} $A:= A_{p}\cdots A_2 A_1$ is strictly positive. These IETs are called in the literature  \emph{periodic-type} IETs (see e.g.~\cite{SU:wea}), in analogy with \emph{periodic-type} rotation numbers (quadratic irrationals like the golden mean $(\sqrt{5}-1)/2=[1,1,\dots, 1, \dots]$ which have a periodic continued fraction expansion). By construction they are \emph{self-similar} and one can also show that they arise as Poincar{\'e} section of foliations which are fixed by a \emph{pseudo-Anosov} surface diffeomorphism. Notice that $d$-IETs of periodic type form a \emph{measure zero set} in $\mathcal{I}_d$ (they are actually countable).  One can show (in view of a Perron-Frobenius argument, e.g.~following \cite{Ve:gau}) that periodic-type IETs are always \emph{uniquely ergodic} with respect to the Lebesgue measure. 

Periodic-type IETs are often the very first type of IETs used to construct \emph{explicit examples}: see e.\ g.\ the explicit examples of weak mixing periodic-type IETs in \cite{SU:wea} or the explicit examples of Roth-type IETs build in the Appendix of \cite{MMY:coh}.   %{\color{blue}More to ADD?}  
%, or \cite{Ul:the} for absence of mixing for locally Hamiltonian flows in $\mathcal{U}_{min}$, see \S~\ref{sec:locHam}). 
On the other hand, among periodic-type IETs one can also find examples with exceptional behavior. %{\color{blue} Refs to ADD?}
 A further request, that is used to guarantee that a periodic-type $T$ display features similar to those of typical (in the measure theoretical sense) IETs, is that $T$ is of \emph{hyperbolic periodic-type}: this means that the periodic matrix $A$ has $g$ eigenvalues of modulus greater than $1$, where $g$ is the genus of the surface of which $T$ is a Poincar{\'e} section. Notice that  $g$ is the largest possible number of such eigenvalues, as it can be shown by either geometric or combinatorial arguments (in particular exploiting the symplectic features of the cocycle matrices, which come from their interpretation as action of renormalization on the \emph{relative} homology $H_1(S, Fix(\varphi_\R), \mathbb{R})$, one can show that $A$ has also $g$ eigenvalues of modulus less than $1$, while the transpose $A^T$ acts as a permutation on a subspace of dimension $k:= d-2g$ which gives rise to a $k$-dimensional central space). 
%We remark that periodic-type IETs form a  \emph{countable} set.

\subsubsection*{Bounded-type IETs equivalent characterizations.}
Periodic-type IETs are a special case of so called {bounded-type} IETs: we say that a (Keane) IET $T$ is of \emph{bounded-type} if the matrices of the \emph{positive acceleration} $P_k:=A(n_k,n_{k+1})$ are uniformly bounded, i.e.~there exists a constant $M>0$ such that $\Vert P_k\Vert\leq M$ for every $k\in\mathbb{N}$. From this point of view, bounded-type IETs can be seen as a generalization of \emph{bounded-type} %{\color{blue}TO CHECK (also called in the literature \emph{constant-type})} 
 rotation numbers (which, recalling   \S~\ref{sec:g1}, are $\alpha=[a_0,a_1,\dots, a_n,\dots ]$ such that for some $M>0$ we have $|a_n|\leq M$). It turns out that this renormalization-based definition characterizes a natural class of IETs (and corresponding surfaces) from the combinatorial and geometric point of view: bounded-type IETs are \emph{linearly recurrent} (i.e.\ satisfy an important notion of low complexity in word-combinatorics) and  surfaces which have a bounded-type IET as a section give rise to \emph{bounded} Teichmueller geodesics in the moduli space of translation surfaces (see e.g.~\cite{HMU:Lag} for the proof of the equivalences). These natural characterizations show once more how the \emph{positive} acceleration (and not simply Zorich acceleration) is the good one to use in this setting (see also \cite{KM:bt} where it is shown that asking that Zorich matrices are bounded leads to a different, strictly larger class). 

Furthermore, from the point of view of renormalization, the uniform bounds on the norm of the matrices $P_k$ imply that  the partitions into Rohlin towers produced by Rauzy-Veech renormalization are all \emph{balanced} (see \S~\ref{sec:renorm}).
From a purely dynamics perspective, 
% from the point of view of renormalization, bounded-type IETs are such that and 
 the orbits of a bounded-type IET are \emph{well-spaced}: there are uniform constants $c, C>0$ such that, for any point $x$ and any $n$, the \emph{gaps} (i.e.~the distances between closest point) of the orbit $\{ T^i x, 0\leq i< n\}$ are all comparable to  $n$, i.e.~are bounded below by $c/n$ and above by $C/n$. Yet another characterization is in terms of orbits of discontinuities:  if $\delta_n(T)$ denotes the  smallest length of a continuity interval for $T_n$, $\liminf_{n\in \mathbb{N}} n\delta_n(T)>0$, see \cite{HMU:Lag} and the reference therein.  
 
Several results in the literature were proved first assuming bounded-type (for example  absence of mixing for flows in $\mathcal{U}_{min}$, see    \cite{Ul:the}, preceding \cite{Ul:abs}) 
% or the Switchable Ratner property and multiple mixing for flows in $\mathcal{U}_{\neg min}$ (see the work \cite{KK:Rat} by Kanigowski and Ku{\l}aga-Przymus, preceeding \cite{KKU:mul}). S
and some properties are currently known only  under the assumption of being bounded-type, for example \emph{absence of partial rigidity} and \emph{mild-mixing} (see \cite{Ku:sel} and \cite{KK:Rat} respectively) for flows in $\mathcal{U}_{min}$ (it is possible, but an open question, that these two properties fail without assuming that a Poincar{\'e} section is of bounded-type), or ergodicity of typical skew-product extensions of IETs by piecewise constant cocycles (see \cite{CR:erg}). %{\color{blue}ADD MORE? MORE Examples?}

\subsubsection*{Bounded-type uniform contraction and deviations estimates}
One of the way in which the bounded-type assumption can be exploited is the following. It is well known that iterates of a \emph{positive} $d\times d$ matrix $A>0$ act on the positive cone $\mathbb{R}_+^d$ as a \emph{strict} contraction (e.g. with respect to the Hilbert projective metric): this is the phenomenon behind the proof of Perron-Frobenius theorem, that shows that $A$ has a unique (positive) eigenvector with maximal eigenvalue. More generally, the projective action of any matrix $A_i$ with $\Vert A_i\Vert \leq M$ has a contraction rate which depends on $M$ only; this, in view of the 
connection between the entries of the cocycle products $A^n:= A_n \dots A_1$ and (special) Birkhoff sums  (see \S~\ref{sec:renorm}), can be used, given a bounded-type IET, to prove unique ergodicity and to give uniform estimates on the rate of convergence of ergodic averages: one can  for example show  that there is a uniform constant (which can be taken to be $1$) and a uniform exponent $\gamma_M$ such that, for any bounded-type IET with $\Vert P_k\Vert\leq M$ and any mean zero (piecewise) smooth $f:I\to \R$, $| S_n f (x) | \leq n^{\gamma_M}$  for all $x\in I$ (see the Appendix of \cite{CR:erg}). 

\subsubsection*{The role of bounded-type in the study of Lagrange spectra}
Periodic-type and bounded-type rotation numbers play a central role in the study of the \emph{Lagrange spectrum} $\mathcal{L}\subset \R \cup \{+\infty\}$, a classical object in both number theory and dynamics (see for example~\cite{HMU:Lag} or \cite{Ma:Lag}  and the reference therein). It is defined as the set $\mathcal{L}:= \{ L(\alpha),\ \alpha \in \R\}$ where $L(\alpha):= \limsup_{q,p\to \infty} 1/{q|q\alpha -p|}$; one can show that $L(\alpha)<\infty$ exactly when $\alpha$ is of bounded-type, in which case  $L(\alpha)^{-1}$ provides the smallest constant such that  $|\alpha - p/q| < L(\alpha)^{-1} /q^2$ has infinitely many integer solutions  $p,q\in\mathbb{Z}$, $q\neq 0$ (and it has  also an interpretation in terms of depths of excursions into the cusp of hyperbolic geodesics on the modular surface). %and it also has an interpretation
%This leads to the introduction of the (classical) \emph{Lagrange spectrum} 
%$\cL \subset \RRbar:=\RR \cup \{+\infty\}$  collection
% of values
 %$\{ L(\alpha), \alpha \in \RR\}\subset \RRbar:=\RR \cup \{+\infty\}$ of values $L(\alpha)\in \RRbar $ given by 
%For a given $\alpha\in\RR$, let   be such that 
%\[
%L(\alpha):= \sup \{ k : |\alpha - p/q| < 1/k q^2 \ \text{for infinitely many } q\in \mathbb{N},\  p \in \mathbb{Z} \} \limsup_{q,p\to +\infty} \frac{1}{q|q\alpha -p|},\quad \alpha \in\RR \right\} .
%\]
% for any $c>L$, we have $|\alpha - p/q| > 1/cq^2$ for all $p$ and $q$ large enough and $L(\alpha)$ is minimal with respect to this property (for example, $L(\frac{1+\sqrt{5}}{2}) = \sqrt{5}$).  
%Then  $\cL$  is the collection of values $\{ L(\alpha), \alpha \in \RR\}$. %if and only if there exists $\alpha\in\RR$ such that $L(\alpha)=L$.  
%Equivalently, one can also write (see for example~\cite{MatheusCIRM})
%\begin{equation}\label{eq:classicdefLagrange}
%	\cL = \left\{ L(\alpha)= \limsup_{q,p\to +\infty} \frac{1}{q|q\alpha -p|},\quad \alpha \in\RR %\right\} \subset \RRbar=\RR \cup \{+\infty\}.%
%\end{equation}
%One can see that for almost every $\alpha$ one has $L(\alpha)=\infty$, but $L(\alpha)<\infty$ for a set of full Hausdorff dimension, which consists exactly of so called \emph{badly approximable} (or bounded type) numbers.  \
Among the many geometric and dynamical extensions of the notion of Lagrange spectrum (see some of the references in \cite{HMU:Lag}), a natural generalization to higher genus leads to Lagrange spectra of IETs and translation surfaces, which we introduced in joint work with Hubert and Marchese in \cite{HMU:Lag}. The finite values of these spectra are achieved exactly by bounded-type IETs and can be computed using renormalization. We show furthermore in \cite{HMU:Lag} that these spectra can be obtained as the \emph{closure} of the values achieved by periodic-type IETs.

\subsection{Roth-like conditions and type}\label{sec:Roth} The \emph{Roth-type} condition, to the best of our knowledge, was historically the first full measure 'arithmetic' condition to be defined and exploited in higher genus.
\subsubsection*{Roth-type condition}
 In the seminal paper \cite{MMY:coh}, Marmi, Moussa and Yoccoz show first of all that (a predecessor of) the positive acceleration of Rauzy-Veech induction (refer to \S~\ref{sec:renorm}) is well defined for all Keane IETs and use this acceleration to define the Roth-type condition and prove that it has full measure; they then  show that this condition is sufficient to solve the cohomological equation after removing obstructions (see  \S~\ref{sec:linhigherg}). Since bounded-type IETs have measure zero, to describe a full measure set of IETs one needs to allow the norms $\Vert P_k\Vert$ of the matrices $(P_k)_k$ of the positive acceleration to grow. Marmi, Moussa and Yoccoz show in \cite{MMY:coh} that, for almost every  $d$-IETs in $\mathcal{I}_d$, the matrices $(P_k)_k$  grow sub-polynomially, i.e.~for any $\epsilon>0$ there exists $C_\epsilon>0$ such that
\begin{equation}\label{eq:Rothgrowth}
\Vert P_k\Vert\leq C_\epsilon \Vert Q_k\Vert^\epsilon,\qquad \text{where}\ Q_k:= P_{k-1}\dots P_{0}.
\end{equation}
This condition should be seen as a higher genus generalization of the classical Roth-type condition, see \S~\ref{sec:g1}. A $d$-IETs is called \emph{Roth-type} if it satisfies \eqref{eq:Rothgrowth} (which is equivalent to condition (a) in  \cite{MMY:coh}, see \cite{MY:Hol}), and two additional conditions,  which concern the contraction properties of the cocycle (condition (b) in \cite{MMY:coh} impose that the operators $S(k)$ act as contractions on mean zero functions and guarantees unique ergodicity and the existence of a \emph{spectral gap}, 
%requests (since it in particular implies unique ergodicity and the existence of gap between the first and second Lyapunov exponent) is that the products $Q_k$ contract uniformly,
 while the last one, condition (c) or \emph{coherence}, concerns the contraction rate of the stable space and its quotient space). The presence of additional requests that concern not only the growth of the matrices but also their hyperbolicity properties seems to be an important and new feature of several Diophantine-like conditions in higher genus, see \S~\ref{sec:OsDC}. While the proof that the last two conditions are satisfied by  almost every IET is a simple consequence of Forni's work \cite{Fo:dev} and Oseledets theorem (which can be applied in view of the work by Zorich \cite{Zo:gau}), the proof that the growth condition \eqref{eq:Rothgrowth} is typical takes a large part of \cite{MMY:coh}; a simpler proof 
%which mimics the proof that the classical Roth-type has full measure using the Gauss map and continued fraction, 
can be now deduced (as explained in \cite{MMM:Yo}) from a later result by Avila-Gouezel-Yoccoz \cite{AGY:exp}. 

\subsubsection*{Variations of the Roth-type condition.}
In a similar way in which one can refine the  periodic-type condition by defining \emph{hyperbolic} periodic-type, one may further request, given a Roth-type IET $T$, that the stable space, i.e.~the space 
$
\Gamma_s(T)$ of vectors $ v\in \R^d$ such that $ A^n v \to 0$ exponentially as $n $ grows (which, in the case of a periodic-type IET with period matrix $A$, is generated by the eigenvectors of the eigenvalues of $A$ which have modulus greater than $1$) has maximal dimension, namely $g$. The condition that one gets was called \emph{restricted Roth-type} in \cite{MMY:lin}; it has full measure in view of \cite{Fo:dev}  and was used to study the structure and codimension of local $\mathcal{C}^r$ conjugacy class of a (G)IET for $r>2$. In the joint work \cite{MUY:onR} with Marmi and Yoccoz, we introduced a further weakening of the (restricted) Roth-type condition, the \emph{absolute} (restricted) Roth type condition, expressed only in terms of the cocycle action on a $2g$ dimensional subspace which can be identified with the \emph{absolute} homology $H_1(S,\mathbb{R})$ of the surface $S$ of which $T$ is section (in contrast, the original condition involves the whole cocycle, which describes the action on the \emph{relative} homology $H_1(S,Fix(\varphi_\R),\mathbb{R})$). Exploiting \cite{EC:Ose}, one can  also show that this absolute (restricted) Roth type condition holds on every translation surface for almost every direction (see \cite{EC:Ose} and \cite{MUY:onR}). A generalization of the \emph{restricted} Roth type condition, the \emph{quasi-Roth type} condition, was introduced in \cite{FMM:coh} to extend the results of \cite{MMY:coh} and \cite{MMY:lin} to Poincar{\'e} maps of surfaces for which the stable space has dimension less than $g$  (see \cite{FMM:coh} for details). 
 Let us also mention that a Roth-type condition can also be imposed on the \emph{backward rotation number} (of a translation flow), requesting a growth rate similar to \eqref{eq:Rothgrowth} for the \emph{dual} cocycle. The corresponding \emph{dual Roth-type condition} was used in \cite{MUY:onR} to study the asymptotic oscillations of the error term in \eqref{eq:expansion} (which we  describe in terms of a \emph{distributional} cocycle- or \emph{distributional limit shape}, see \cite{MUY:onR} for details).
% and proved that it is satisfied by a full measure set of $d$-IETs in $\mathcal{I}_d$ and showed 

\subsubsection*{Type and recurrence for IETs}
It is not surprising that Diophantine-like conditions can also be used to study \emph{recurrence} questions. While for rotations these reduce to Diophantine properties in the classical arithmetic sense (namely how well a number can be approximated by rationals), given an IET $T$ one can study either how frequently the successive iterates $(T_n(x))_{n\in \N}$ return close to $x$ (see e.g.~\cite{BC:Dio}), or how close the iterates of a discontinuity come to other discontinuities, see e.g.~%the generalization of Khintchin's theorem by Marchese 
\cite{Ma:Khi}. The (Diophantine) \emph{type} $\eta$ of a rotation $R_\alpha$ is defined to be $\eta:=\sup \{ \beta$ s.~t.~$\liminf_{n\to \infty} n^\beta \{ n\alpha \}=0\}$. Bounded-type and Roth-type numbers have type $\eta=1$ (while Liouville ones have type $\eta=\infty$). One can show (see \cite{Ki:typ} and \cite{MMY:lin}) that requesting that an IET $T$ is Roth-type is equivalent to asking that  $\sup \, \{ \beta$ s.~t.~$ \liminf n^\beta \delta_n(T)=0\}=1,$
%$$\sup \{ \beta \ \text{such\ that}\  \liminf n^\beta \delta_n(T)=0\}=1,$$ 
where here $\delta_n(T)$ is the minimum spacing between discontinuities of $T_n$.  
%oth problems have been intensively studied, see for example ADD (Marchese and Boshernitzan Chaika). 
It also implies (but without equivalence)  that the first return time $\tau_r(x)$ of $x$ to a ball of radius $r>0$ satisfies the logarithmic law $\lim_{r\to 0}\log\tau_r(x)/\log (1/ r)=1$ for almost every $x\in [0,1]$ (see \cite{Ki:typ}).

\subsection{Controlled growth  Diophantine-like conditions.} \label{sec:growth}
Any \emph{balanced acceleration} of Rauzy-Veech induction (as defined in \S~\ref{sec:renorm}), produces, given a typical IET $T$, a  sequence of times $(n_k)_k$ which correspond to occurrences of positive matrices $A_{n_k}$ whose norm $\Vert A_{n_k}\Vert \leq M$ is uniformly bounded
 (these are furthermore return times  to a compact subset $K$ of the parameter space for the natural extension).  As for bounded-type IETs,  occurrences of these positive bounded matrices give very good control of the convergence of (special) Birkhoff sums of characteristic functions $\chi_{I_j}$ (see the end of \S~\ref{sec:bt}). More generally, if  $x_0\in I_{n}^j$ belongs to the inducing interval $I_n$ of a balanced return time $n:=n_k$ and   $q	:=r_{n}^j$ is the height of the corresponding tower, the orbit  $\{ x_0, T(x_0), \cdots , T^{q-1}(x_0)\}$ \emph{along a tower} is  so regularly spaced that 
one can  get good estimates of the Birkhoff sums  $S_{q} f (x_0)$ 
%if $q_{n_k}^j$ are return times of points $x\in I^n_j$, 
also for  other classes of observables $f$. In order to estimate Birkhoff sums $S_n f (x)$ for other times $n\in \mathbb{N}$ and points $x\in[0,1]$, one can then \emph{interpolate} these estimates by using the decomposition \eqref{decomposition} into special Birkhoff sums. It is clear now that for this interpolation to provide good estimates for any time $n\in \N$, one needs to impose that the balanced times $(n_k)_k$ are sufficiently frequent so that $\Vert A(n_k, n_{k+1})\Vert$ grows in a controlled way. Notice that by balance the tower heights $r_{n_k}^j$ for $1\leq j\leq d$ are all comparable and if we set $q_n:=\max_j r_n^j$, the norm  $\Vert A(n_k, n_{k+1})\Vert$ is proportional to $q_{n_{k+1}}/q_{n_k}$.
%the ratios $q^{n_{k+1}}_j /q^{n_k}_i$ between successive balanced tower heights are controlled. 

\subsubsection*{Mixing Diophantine Condition.} 
The main requirement of the Mixing Diophantine condition introduced in \cite{Ul:mix} is that there exists a (good) positive acceleration and $C>0$ s.t.\  
\begin{equation}\label{eq:MDC}
\Vert A(n_k, n_{k+1})\Vert \leq C k^\tau, \quad \forall\ k\in \mathbb{N}, \quad \text{for\ some}\ 1<\tau<2.
\end{equation}
This condition should be seen as a higher genus generalization of the Khanin-Sinai condition $|a_k|\leq C k^{\tau}$ for mixing of Arnold flows, see \S~\ref{sec:g1}. The proof that it is satisfied by a full measure set of IETs follows from a Borel-Cantelli argument analogous to the one that can be used in genus one, but % while in the classical case it uses that $\int_{[0,1]} a_0^{\nu}(x)\d \mu_\G $ is finite for any $\nu<1$,
 the input in higher genus are the highly non-trivial integrability estimates for balanced accelerations proved by Avila, Gouezel and Yoccoz (which the authors proved to show in \cite{AGY:exp} that the Teichm{\"u}ller geodesic flow is exponentially mixing): it is proved in   \cite{AGY:exp} that for any $0<\nu <1$, there exists a suitable compact set $K$ such that $\int_K \Vert A_K \Vert^\nu \mathrm{d}\mu $ is finite (where $A_K$ is the accelerated cocycle and $\mu$ the Zorich measure). 

In order to prove mixing of (minimal components of) locally Hamiltonian flows in $\mathcal{U}_{\neg min}$ (i.e.~Theorem~\ref{thm:mixing}), one needs good quantitative estimates on \emph{shearing}: these are given by estimates of Birkhoff sums $S_n f$ over an IET which arise as Poincar{\'e} map, for a particular observable $f$ (namely, $f$ is taken to be the derivative of the roof function in the special flow representation of $\varphi_\R$), which turns out not to be in $L^1$ (indeed the function $f$ has singularities of type $1/x$, which are not integrable).  When  $n=n_k$ is a balanced time, on can control the corresponding special Birkhoff sums $S(n_k) f$ and show that each Birkhoff sum along a tower   $S_{q}f(x)$ where $q=q_{n_k}^j$ and $x\in I_{n_k}^j$ can be controlled after  removing the \emph{closest point} contribution that, in this case, is simply $1/x$. One can indeed show that the  
\emph{trimmed} Birkhoff sum $S_{q} f(x)-1/x$ is  asymptotical to $C q \log q$. The Mixing Diophantine Condition allows to \emph{interpolate} these estimates and show that, also for any other $n\in\mathbb{N}$,  $S_n f(x)$ grows asymptotically $C n\log n$ for all points $x$ \emph{with the exception} of points which belong to a set $\Sigma_n\subset [0,1]$ of measure going to zero. The set $\Sigma_n$ of points which needs to be removed to get the desired control contains points whose orbits may be \emph{resonant}, in the sense that it may contain a close-to-arithmetic progression  near one of the singularities of $f$, with step $q_{n_k}/q_{n_{k+1}}$ (which can be a very small step if $q_{n_{k+1}}$ is much larger than $q_{n_k}$).

\subsubsection*{Ratner Diophantine condition.}
In order to prove that (minimal components of) locally Hamiltonian flows in $\mathcal{U}_{\neg min}$ have the \emph{Switchable Ratner property} (e.g.~Theorem~\ref{thm:Ratner}, see~\S~\ref{sec:higherg}), one needs more delicate quantitative shearing estimates. Such estimates are proven assuming first of all the Mixing Diophantine Condition, but the MDC is  not sufficient. While mixing is an asymptotic condition and therefore it is sufficient, for all large $n$, to prove estimates for the Birkhoff sums $S_n f(x)$ (introduced in the previous subsection) on sets of measure tending to $1$ (and hence one can remove a  set $\Sigma_n$ whose measure goes to zero), the (Switchable) Ratner Property requires estimates on arbitrarily large sets of initial points, for \emph{all} large times $n\geq n_0$. 
If the series   $\sum_{n\in \mathbb{N}} Leb(\Sigma_n)$ were finite, the tail sets of the form $\bigcup_{n\geq n_0} \Sigma_n$ would have arbitrarily small measures and thus one could throw away these unions for $n_0$ large. Unfortunately, one can check that the measures  $\big( Leb (\Sigma_n)\big)_{n\in \N}$ are \emph{not} summable. Instead, we consider a subset $K\subset \mathbb{N}$ such that  $\sum_{n\notin K} Leb(\Sigma_n)<+\infty$ and exploit the additional freedom given by the \emph{switchable} Ratner condition to deal with points $x\in \Sigma_n$ when $n\in K$. This requires the introduction of a suitable Diophantine-like condition.
 
We say that an IET $T$ satisfies the \emph{Ratner Diophantine condition} (RDC) if $T$ satisfies the Mixing DC along the sequence $(n_k)_{k\in\mathbb{N}}$ of balanced induction times and if  
%for some  $\nu>1$ and $\overline{l}\in \mathbb{N}$ and 
%there exists $1< \tau < 2$ such that $T \in \MDCo{\tau}$ satisfies the mixing Diophantine condition along a sequence $\{ n_\ell\}_{\ell\in \mathbb{N}}$ and  
there exists $0< \xi, \eta<1$ such that, if  $B_k:=A(n_k,n_{k+1})$  are the matrices of the accelerated cocycle and $q_{k} :=\max_j r_{n_k}^j$ the maximum height of the corresponding towers, we have
%(\emph{Ratner DC}) \emph{with integrability power} $\tau\in(1,2)$ and exponents $(\xi,\eta)$ 
%\todoA[inline]{We can put it in the parathesis.}
%if it satisfies the mixing Diophantine condition with power $\tau$ along a sequence $\{ n_\ell\}_{\ell\in \mathbb{N}}$ of induction times and for $\xi,\eta<1$ 
\begin{equation} \label{eq:RatnerDC}
\sum_{k\notin K } {1}/{{(\log q_k)}^\eta } < + \infty, \qquad \text{where}\  
K:=\{ k\in\N \ \text{s.t.\ } \|B_k\|\leq  k^\xi  \}.
% L=\overline{\ell}(1+ \left[ \log_d (2 \nu^2) \right]).
\end{equation}
The assumption \eqref{eq:RatnerDC} guarantees in particular the summability of $\sum_{n\notin K} Leb(\Sigma_n)$, so that tail sets  of this series \emph{can} be removed. When $k\in K$, using that $n_k$ is a  balanced time  and  $q_{k}/q_{k-1}\leq \|B_k\|$ is not too large, one can show that
an arbitrarily large set of points $x$ do not get close of order $c/q_{k-1}$ to a singularity  \emph{twice}  in time of order $q_{k}$, so either going \emph{forward} or \emph{backward} in time one can avoid getting  $O(q_{k}^{-1})$ close to  singularities. 
% and control $S_{} f(x)$. 
%if an orbit of a point of length $q_\ell$ gets too close to a singularity in the future (where too close is of order $q_{\ell+L}$ for a fixed $L$), then it did not come that close to a singularity in the past
% one can show that if Birkhoff sums are not controlled in the future, they are controlled in the past (Lemma \ref{prty}). 
This suffices to  provide the control of $S_n f (x)$ (and therefore of \emph{shearing})  required by the switchable Ratner property for all times. 

%When $n\in K$, we control the behavior of orbits along towers (and hence the corresponding Birkhoff sums) exploiting the versatility  of the \emph{switchable} Ratner property. 
%One is then left to estimate the times $\ell \in \widetilde{K}_T$.
 %This is where one exploits the versatility of the \emph{switchable} Ratner property, according to which, if the desired quantification of parabolic divergence does not hold for  \emph{forward} Birkhoff sums (see (i) in Definition \ref{def:SWR}), one can \emph{switch} the direction of time, i.e. prove quantiative divergence estimates on \emph{backward} Birkhoff sums (see (ii) in Definition \ref{def:SWR}).  

  Notice that if an IET $T$ is of bounded type (so $\Vert B_k\Vert$ are bounded) then the RDC is automatically satisfied (since the complement of $ K$ in $\N$ is finite and therefore the series is a sum of finitely many terms). The Ratner DC imposes that the times $k$ for which $\Vert B_k \Vert$ is large are not too frequent: in a sense if an IET satisfies the RDC, it behaves like an IET of bounded type modulo some error with small density (as a subset of $\mathbb{N}$), but this relaxation allows the property to hold for almost every IETs: we prove indeed in  \cite{KKU:mul}  that, for suitable choices of $\xi$ and $\eta$ the RDC is satisfied by a full measure set of IETs. Formally (when using the suitable acceleration), the assumption \eqref{eq:RatnerDC} looks like the Diophantine Condition for rotations introduced by Kanigowski and Fayad in \cite{FK:mul}, see \eqref{alphaRatner}. The proof of full measure of the RDC is modeled on the proof of full measure of the arithmetic condition \eqref{alphaRatner}, with the role of the Gauss map played by the renormalization operator in parameter space corresponding to  the balanced acceleration. Key ingredients to make this proof work are once more the integrability estimates from Avila-Gouezel-Yoccoz \cite{AGY:exp}, as well as a \emph{quasi-Bernoulli} property of the balanced acceleration, see \cite{KKU:mul} for details.

\subsubsection*{Backward growth condition for absence of mixing}
The Diophantine-like condition to prove absence of mixing of typical locally Hamiltonian flows in $\mathcal{U}_{min}$ (see \S~\ref{sec:locHam}) is not explicitly stated in \cite{Ul:abs}, but, from the proof, one can see that one needs the existence of a suitable acceleration of the balanced acceleration, whose matrices will be denoted by $(B_k)_{k\in\N}$, of a subsequence $(k_l)_{l}$ and of a constant $M>0 $ such that  
\begin{equation}\label{eq:series}
 \sum_{k=0}^{k_l}\frac{ \Vert {B_{k}} \Vert } { \nu^{k_l-k}} = \sum_{j=0}^{k_l} \frac{\Vert {B_{k_l-j}}\Vert } { \nu^{j}} \leq 
 M <+\infty, \qquad \text{for\ all}\ k\in \mathbb{N}, 
\end{equation}
where $\nu$ is some constant with $\nu>1$. This type of condition has two interesting features: it requires a \emph{backward} control of the growth of the matrices of an accelerated cocycle, which has to happen \emph{infinitely often}. Indeed, for the series \eqref{eq:series} to converge and be uniformely bounded by $M$, one needs to ask that the norms $\Vert B_{k} \Vert$ when $k$ belongs to the sequence $(k_l)_{l\in\N}$ are uniformly bounded; furthermore, it is sufficient to then impose that, going backward in time, they grow less fast than the denominator, namely that $\Vert B_{k_l-j} \Vert \leq C e^{\delta j}$ for $0\leq j\leq k_l$ where $\delta$ is chosen so that $e^\delta<\nu$. These conditions can be shown to be of full measure by exploiting Oseledets integrability (for the \emph{dual} cocycle).  

These type of backward conditions seem to appear naturally when one wants to provide good control of the deviations of the points in a finite segment $\{ x, T (x), \dots, T_N(x)\}$  of an IET orbit from an arithmetic progression: one would like to show for example that, if we relabel the points in the orbit segment so that $0<x_1 <x_2< \dots < x_N<1$, the points $x_i$ display polynomial deviations from an arithmetic progression, i.e.~there exists $C>0$ and $0<\gamma<1$ such that  $|x_i - i/N |\leq  C(i/N)^\gamma$.  
These estimates (which are used in \cite{Ul:the, Ul:abs} to show, through a cancellations mechanism, that there is a subsequence of times with no shearing and, as a consequence, that mixing fails) can be proved for all times for bounded type IETs (see \cite{Ul:the}), but, for typical IETs, even for orbits along a balanced tower of some renormalization level $n_{k_0}$, it may not be possible to choose a constant $C$ uniformly on $i$. Heuristically, the reason for this is that to estimate the location of $x_i$ one can use a \emph{spatial decomposition} of the interval $[0,x_i]$ into floors of renormalization towers  which involves the entries of \emph{backward} cocycle matrices (a decomposition similar to the one in  \eqref{decomposition}, but with the role of time now played by space; geometrically this can also be interpreted as swapping the role of the horizontal and vertical flows on a translation surface). The presence of an exceptionally large $\Vert A_k\Vert$, even if $k$ is much smaller than $k_0$, can still spoil the deviations control, since it may correspond in the spatial decomposition to a \emph{clustering} of points, close to an arithmetic progression of a very small step. % and this may spoil the desired control.   %each $x_i$ belongs indeed to a the base interval $I_{n_k}$ of \emph{level} $n_k$ of the acceleration towers 
 
We point out that phenomenona of similar nature, where the whole backward history of the continued fraction entries matters to control orbits, appears also in genus one,  in the theory of circle diffeomorphisms. In the paper \cite{KT:Her}, in which Herman's theory of linearization (see \S~\ref{sec:g1}) is revisited through renormalization following \cite{KS:Her} and optimal results are achieved for low regularity, the arithmetic condition required on $\alpha=[a_0,a_1,\dots]$ is a condition involving a series, namely the finiteness of $\sum_{n=0}^\infty a_{n+1} \big(\sum_{i=0}^n \frac{l_n}{l_{n-i}} ({l_{n-i-1}})^\eta\big)$, where $l_n:= |q_n\alpha-p_n|$ and $0<\eta<1$. 
This condition is then also in \cite{KS:Her} to control the \emph{spatial decomposition} of orbit segments. It would be interesting to know if the analogy, which  at this level is only formal and on the \emph{nature} of the conditions, hides a more profound similarity.

\subsection{Effective Oseledets Diophantine-like conditions.} \label{sec:OsDC}
To conclude, we briefly describe the Regular and Uniform Diophantine-like conditions (RDC and UDC for short), introduced and used to prove Theorem~\ref{thm:deviations} and Theorem~\ref{rigidityg2} respectively (see \S~\ref{sec:locHam} and \S~\ref{sec:linhigherg}). 
 %respectively in \cite{FU:Bir} to control deviations of ergodic averages for locally Hamiltonian flows (see \S~\ref{sec:locHam}, Theorem~\ref{thm:deviations}) and  to prove a priory bounds for GIETs and rigidity in genus two (see \S~\ref{sec:linhigherg}, Theorem~\ref{rigidityg2}). 
	Both these conditions present a novel aspect: not only they impose a controlled \emph{growth} of cocycle matrices of a suitable acceleration (as all the conditions we have seen in \S~\ref{sec:growth}), as well  \emph{hyperbolicity} assumptions (see the \emph{hyperbolic} periodic-type or the \emph{restricted} Roth type condition, in \S~\ref{sec:bt} and \S~\ref{sec:Roth}), but they also impose \emph{quantitative} forms of \emph{hyperbolicity}, %m The latter aspect is controlled 
by %	imposing  Oseledets genericity with 
asking for	\emph{effective} bounds on the convergence rates in the  conclusion of Oseledets theorem, as we now detail.

\subsubsection*{Effective Oseledets control and the UDC}\label{sec:UDC}
Let us say that a sequence of balanced return times $(n_k)_{k\in \N} $ 
%and the corresponding accelerated matrices $B_k:= A(n_k, n_{k+1})$ 
satisfies an {\emph{effective} Oseledets control}  if one can find a sequence of \emph{invariant splittings} $\mathbb{R}^d=E^n_s\oplus E^n_c\oplus E^n_u$, with $dim E^n_s=g$, such that, for some $\theta>0$ and any $k\in \mathbb{N}$% we have 
\begin{align}
 ||A(n_k,n)\vert_{E_s^{n_k}}||_{\infty} \leq C e^{-\theta(n-n_k)} & \qquad \textrm{for\ every} \ n\geq n_k; \label{eq:contr} \\
 ||A(n,n_k)^{-1}\vert_{E_u^{n_k}}||_{\infty}  \leq C e^{-\theta(n_k-n)} & \qquad \textrm{for\ every}\  0\leq  n\leq n_k.\label{eq:exp}
\end{align}
Thus, the cocycle contracts the stable space $E_s^{n_k}$ in the future and the unstable space $E_u^{n_k}$ in the past with a \emph{uniform} rate $\theta$ and a uniform constant $C$. These times can be produced for example considering returns to a set (for the natural extension) where the conclusion of Oseledets theorem (for the cocycle and its inverse) can be made uniform. 
% at all times $(n_{k_l})_l$. 
%The Diophantine-type condition that we will use in the main theorems is the following.
%\begin{definition}[\customlabel{UDC}{UDC}]\label{def:UDC}
An IET satisfies the \emph{Uniform Diophantine Condition} (UDC) if there exists balanced times $(n_k)_k$ with effective Oseledets control 
 and furthermore, for every $\epsilon >0$ there exists $C,c>0, \lambda>0$ and  a subsequence  $(k_l)_{l\in \N}$ which is \emph{linearly growing}  (i.e.~such that $\liminf_{l\to \infty} k_l/l>0$) for which
%acceleration for which $T$ is Oseledets generic with effective Oseledets control and,  
%furthermore, there exists constants $C >0$ and $\tau_1>\tau_0>0$ and, for any $\epsilon>0$, %$C_\epsilon>0$ such that %the matrices $B_k$  
%\begin{align}
%\tag{$O$} \quad  & T\  \text{is\ \emph{Oseledets\ generic}}, \text{i.e.}\ there exsits \label{def:O} \ \text{admits} \ \\ &  \text{an %Oseledets}\ \text{splitting of hyperbolic type, as in \S~\ref{sec:Oseledets};}\nonumber%
%\end{align}
% $T$ is hyperbolic (in the sense of Definition~\ref{def:Rokhlinbalanced}) for the accelerated cocycle $A_Y$
%and the product matrices $B^k:=A(0,n_k)$ satisfy
 %there exist  constants $\kappa>1$ and $\lambda_1>\ldots>\lambda_g>0$ such that
\begin{align}
%&\forall \ \epsilon>0, \ \exists \ C_\epsilon>0, \ s.t.\quad 
& \|A(n_k,n_{k_l})\|\leq C_\epsilon e^{\epsilon |k -k_l|} \text{ for all } k\geq 0\text{ and } l\geq 0;\label{eq:subexpUDC}
%C_\epsilon e^{\epsilon k}
%\|B_k\|\leq 
\\
& c e^{\lambda k}\leq\|A(0,n_k)\|\leq Ce^{(\lambda+\epsilon) k} \text{ for all }\  k\geq 0.\label{eq:growthUDC}
\end{align}
One can show  that assuming that $T$ satisfies the UDC implies in particular that $T$ is of (restricted) Roth-type (see \cite{FU:Bir}); on the other hand, \eqref{eq:contr} and \eqref{eq:exp} are assumptions of a new nature, and furthermore \eqref{eq:growthUDC} clearly excludes IETs of bounded type;  thus this is a more restrictive Diophantine-like condition, although still full measure (see \cite{FU:Bir}).

\subsubsection*{The RDC and conditions on Diophantine series}
In the Regular Diophantine Condition (used to study rigidity of GIETs in \cite{GU:rig} and in particular to prove Theorem~\ref{rigidityg2}) we assume that $T$ is Oseledets generic and require the existence of a special sequence of balanced times $(n_k)_k$ such that the two following \emph{forward} and \emph{backward} series (involving the accelerated matrices $B_k:=A(n_k,n_{k+1})$, their products $B(k,l):= B_{l-1} \cdots B_k $, as well as the projections $\Pi_s^k $ and $ \Pi_u^k$ to $E_s^{n_k}$ and $ E_u^{n_k}$ respectively) 
%is of similar nature: we also require the existence of a sequence balanced times $(n_k)_k$  
%(corresponding to a double occurrence of a bounded positive matrix) 
%with effective Oseledets control %as given by \eqref{eq:contr} and \eqref{eq:exp} 
%and of a linearly growing subsequence $(k_l)_l$ which satisfy \eqref{eq:subexpRDC}. In addition, along this subsequence, we also require a uniform lower bound on the \emph{angles} between the subspaces $E^n_s, E^n_u$ and $E^n_c$ of the splitting and, furthermore, the following key condition, that the two following  series 
%(that we call \emph{forward} and \emph{backward} series respectively) 
are uniformly bounded by some constant $M>0$ along a linearly growing subsequence $(k_l)_{l\in \N}$, namely, for every $l\in \N$,
%as well as effective control of the angles 
%$$ |\angle (E^{(n)}_x, E^{(n)}_y)|\geq c_2 \ e^{-\epsilon|n-n_k|}, & \qquad \textrm{for\ all}\ n\in\mathbb{Z}, \textrm{distinct}\, x,y\in\{s,c,u\};$$
%{Global aspect}: the exists a uniform $C>0$ such that for all $k$
\begin{equation} \label{eq:bfseries} \sum_{k=1}^{k_k}{ ||B (k,k_l)_{| E_s^{n_k}}|| \,|| \Pi_{s}^{k}|| \,||B_k||} \leq M ,  \quad  \sum_{k=k_l +1}^{\infty}{||B(k_l,k)^{-1}_{| {E}_u^{k}}|| \,||\Pi_{u}^{k}|| \,||B_k||} \leq M .
\end{equation}
 We also require a uniform lower bound on the \emph{angles} between the subspaces $E^n_s, E^n_u$ and $E^n_c$ of the splitting along the subsequence $(n_{k_l})_l$ and sub-exponential growth of $B(k_l,k_{l+1})$.  The convergence of these series can be proved assuming that the sequence $(n_k)_k$ provides effective Oseledets control; the subsequence $(k_l)_l$ is then selected so that the uniform upper bound holds. Also the  UDC can be used to prove the convergence and uniform boundedness (along a linearly growing subsequence) of some series of similar (although simpler) nature (that we call \emph{Diophantine series}, see \cite{FU:Bir} for details). Notice also the similarity between the backward series in \eqref{eq:bfseries} and the series \eqref{eq:series} used to prove absence of mixing, even though the latter involves only the norm of the matrices and not their hyperbolic properties. 

Examples of arithmetic conditions on classical rotation numbers which do not depend only on the asymptotic behavior of the continued fraction entries (as \emph{Diophantine} or \emph{Roth-type} conditions) but instead  depend on values or finiteness of series involving  continued fraction entries include the \emph{Brjuno}-condition (see e.g.~\cite{Yo:pet}) and the \emph{Perez-Marco} condition \cite{PM:sur}. Conditions which require recurrence to a set of rotation numbers with this type of control in the theory of circle diffeos seem to appear in global rigidity results, see for example the Condition $(H)$ defined by Yoccoz (see \cite{Yo:CIME}).  
%{\color{red} Add ref to CIME lecture notes, only place where apparently it is published according to Marmi} 
%as well as the optimal condition for the conjugacy problem in finite smoothness described  by Teplinsky and Khanin in \cite{KT:Her}. 

\subsubsection*{Final remarks and questions} We saw that advancements in our understanding of both chaotic properties and linearization and rigidity questions in the context of surface flows in higher genus depend crucially on sometimes delicate Diophantine-like conditions, imposed to control the renormalization dynamics. While some of these resembles the classical counterparts, others are of new nature, in particular involving hyperbolicity features which become visible only in higher genus. A downside of this new aspect is that conditions that requires Oseledets genericity assumptions are not easily checkable. If there is a way of producing explicit  examples with such properties which are not of periodic type, even within a locus,  remains a challenge.  Since many developments are still quite recent, it is possible that some conditions can be simplified or weakened and still yield the same results; furthermore, the interdependence or inclusions between the various conditions have not been fully investigated.  Finally, even though, all the   conditions we described, with the only exception of bounded-type conditions,  are of full measure, they are likely not  to be the optimal ones required for the results for which they were introduced (we know this for example for the  absence of mixing  condition, in view of \cite{ChW:mix}). Finding optimal conditions for each of these problems is certainly interesting but { probably} very difficult. 

\subsection*{Acknowledgements}
I would like to acknowledge Selim Ghazouani, Giovanni Forni and Stefano Marmi for many interesting discussions, as well as for their useful feedback on an earlier version of this text. A thank goes also to Stefan Kurz for allowing me to use Figure~\ref{renormalization} and to P.\ Berk, K.\ Fr{\c a}czek, D.\ H.\ Kim,  F. Trujillo and T.~Kappeler for reading and commenting this text. 

%\subsection*{Funding}
%This work was partially supported by the Swiss National Science Foundation.
%\end{funding}

%------
% Insert the bibliography.
%------

\bibliographystyle{emss}

%\begin{thebibliography}
\bibliography{biblio_ICM_ulcigrai}
%\end{thebibliography}

\end{document}